\documentclass[onefignum,onetabnum]{siamart190516}
\pdfoutput=1
\usepackage{pdfpages}
\usepackage{mathtools}
\usepackage{autonum}
\usepackage{amssymb}
\usepackage{arydshln} 
\usepackage[utf8]{inputenc}
\PassOptionsToPackage{final}{graphicx}
\usepackage{pgfplots}
\usepackage{ifdraft}
\ifdraft{
\usepackage[colorinlistoftodos]{todonotes}
}{\usepackage[disable]{todonotes}}
\usepackage{dirtytalk}
\usepackage{relsize}
\usepackage{textgreek}

\DeclareMathOperator{\E}{\mathbb{E}}
\DeclareMathOperator{\V}{\mathbb{V}}
\DeclareMathOperator{\Cov}{\mathrm{Cov}}

\renewcommand{\epsilon}{\varepsilon}
\newcommand{\R}{\mathbb{R}}
\newcommand{\CV}{\tilde{v}}
\newcommand{\VD}{\bar{V}}
\newcommand{\PT}{u}

\newcommand{\D}{\mathcal{D}}

\newcommand{\dm}{{\mathrel{\Delta m}}}
\newcommand{\pcf}{p_{c,\Delta t_\ell}}

\newcommand{\pncf}{p_{nc,\Delta t_\ell}}
\newcommand{\pncr}{p_{nc,\Delta t_{\ell-1}}}

\newcommand{\mse}{E}
\newcommand{\dx}{\text{d}x}
\newcommand{\dv}{\text{d}v}

\usetikzlibrary{external}
\tikzexternaldisable

\newcommand\overmat[2]{
  \makebox[0pt][l]{$\smash{\overbrace{\phantom{%
    \begin{matrix}#2\end{matrix}}}^{\text{$#1$}}}$}#2}

\usepackage{soul}

\usepackage{lipsum}
\usepackage{amsfonts}
\usepackage{graphicx}
\usepackage{tikz}
\usepackage{pgfplots}
\pgfplotsset{compat=newest}
\usepackage{epstopdf}
\usepackage{subcaption}
\usepackage{algorithmic}
\ifpdf
  \DeclareGraphicsExtensions{.eps,.pdf,.png,.jpg}
\else
  \DeclareGraphicsExtensions{.eps}
\fi

\newsiamremark{remark}{Remark}
\newsiamremark{hypothesis}{Hypothesis}
\crefname{hypothesis}{Hypothesis}{Hypotheses}
\newsiamthm{claim}{Claim}

\headers{Simulation of Boltzmann-BGK near the diffusive limit with APMLMC}{E.~L{\O}vbak,  and G.~Samaey}

\title{Accelerated simulation of Boltzmann-BGK equations near the diffusive limit with asymptotic-preserving multilevel Monte Carlo\thanks{Submitted to the editors DATE. Some sections of this paper are based on an LNCS proceedings paper referenced as~\cite{Loevbak2020}.
\funding{This work was funded by the Research Foundation - Flanders (FWO) under fellowship number 1SB1919N. The computational resources and services
used in this work were provided by the VSC (Flemish Supercomputer Center), funded by the Research Foundation - Flanders (FWO) and the Flemish
Government - department EWI.}}}

\author{Emil L{\O}vbak\thanks{NUMA Section, Department of Computer Science, KU Leuven, Belgium
  (\email{emil.loevbak@kuleuven.be}, \email{giovanni.samaey@kuleuven.be}).} \and Giovanni Samaey\footnotemark[2]}

\ifpdf
\hypersetup{
  pdftitle={Accelerated simulation of Boltzmann-BGK equations near the diffusive limit with asymptotic-preserving multilevel Monte Carlo},
  pdfauthor={Emil~L{\o}vbak and Giovanni~Samaey}
}
\fi

\definecolor{color0}{rgb}{0.12,0.47,0.7}
\definecolor{color1}{rgb}{0.9,0.5,0.1}
\definecolor{color2}{rgb}{0.3,0.7,0.3}
\definecolor{color3}{rgb}{0.7,0.2,0.8}

\begin{document}

\maketitle

\begin{abstract}
Kinetic equations model the position-velocity distribution of particles subject to transport and collision effects. Under a diffusive scaling, these combined effects converge to a diffusion equation for the position density in the limit of an infinite collision rate. Despite this well-defined limit, numerical simulation is expensive when the collision rate is high but finite, as small time steps are then required. In this work, we present an asymptotic-preserving multilevel Monte Carlo particle scheme that makes use of this diffusive limit to accelerate computations. In this scheme, we first sample the diffusive limiting model to compute a biased initial estimate of a Quantity of Interest, using large time steps. We then perform a limited number of finer simulations with transport and collision dynamics to correct the bias. The efficiency of the multilevel method depends on being able to perform correlated simulations of particles on a hierarchy of discretization levels. We present a method for correlating particle trajectories and present both an analysis and numerical experiments. We demonstrate that our approach significantly reduces the cost of particle simulations in high-collisional regimes, compared with prior work, indicating significant potential for adopting these schemes in various areas of active research.
\end{abstract}

\begin{keywords}
  kinetic equations, multilevel Monte Carlo, particle methods, asymptotic-preserving schemes
\end{keywords}

\begin{AMS}
  65C05, 65C35, 65M75, 35Q20
\end{AMS}

\section{Introduction}

In this work, we consider an ensemble of particles that travel freely until they receive a new velocity due to collisions. Such ensembles are used for modeling physical phenomena in many application domains. An example, which draws our interest, is plasma edge simulation for the design of tokamak fusion reactors, such as performed in the simulation codes B2-EIRENE~\cite{Reiter2005} and DEGAS~\cite{Stotler1994}. Other areas include bacterial chemotaxis~\cite{Rousset2013} and computational fluid dynamics~\cite{Pope1981}, where a large number of software packages exist such as MONACO~\cite{Dietrich1996} and dsmcFoam+~\cite{Scanlon2010,White2018}.

In the cases we consider, the particle ensembles follow a kinetic-transport equation
\begin{equation}
\label{eq:kinetic_abstract}
\partial_t f(x,v,t) + v \cdot \nabla_x f(x,v,t) = Q[f(x,v,t)],
\end{equation}
where $f(x,v,t)$ represents the particle probability density as a function of space $x\in \D_x \subset \R^d$, velocity $v \in \D_v \subset \R^d$ and time $t \in \R^+$, and $Q[\cdot]$ represents a collision operator that causes particles to undergo discontinuous velocity changes. These models are high-dimensional, as they lead to time dependent simulations in a $2d$-dimensional phase-space, e.g., a three dimensional application requires a seven dimensional model. The quantity of interest (QoI) is, however, not the full solution in $2d+1$ dimensions, but rather some lower dimensional function thereof, e.g., the particle position-density $\displaystyle\rho(x,t) = \int_{\D_v} f(x,v,t)\dv$ or momentum $\displaystyle j(x,t) = \int_{\D_v} v f(x,v,t)\dv$ as a function of time.

Aside from being lower-dimensional, these quantities of interest often adhere to slower time-scales than those of the particle dynamics, especially for high collision rates, i.e., the expected number of collisions per second. Having multiple time-scales in the models implies problem stiffness, i.e., explicit methods require small simulation time steps, relative to the simulation horizon, to capture the particle dynamics, resulting in a high simulation cost. When selecting a discretization for models of the form \eqref{eq:kinetic_abstract}, both their high-dimensionality and potential stiffness should be considered.

We consider a diffusively scaled version of \eqref{eq:kinetic_abstract} where the collision operator $Q[\cdot]$ is chosen to be the Bhatnagar-Gross-Krook (BGK) operator~\cite{Bhatnagar1954}, i.e.,
\begin{equation}
\label{eq:kineticdiffusive}
\partial_t f(x,v,t) + \frac{v}{\epsilon} \cdot \nabla_x f(x,v,t) = \frac{1}{\epsilon^2}\left(\mathcal{M}(v)\rho(x,t)-f(x,v,t)\right),
\end{equation}
with $\epsilon$ the Knudsen number, which can intuitively be understood as the inverse of a dimensionless collision rate, and $\mathcal{M}(v)$ the probability density of the post-collisional velocity distribution. The BGK operator affects each particle individually, linearly driving the velocity distribution to a steady-state $\mathcal{M}(v)$ at a rate $\epsilon^{-2}$, rather than modeling particle-particle interactions. As $\epsilon \to 0$, the time-scale separation between the transport and collision behavior in \eqref{eq:kineticdiffusive} becomes infinite. The diffusive scaling, where one simultaneously rescales the transport and collision terms results in convergence to a drift-diffusion equation in this limit~\cite{Dimarco2014}.

As $\epsilon \to 0$, the velocity distribution relaxes to the steady state distribution with density $\mathcal{M}(v)$. In this limit the velocity dimensions of the phase-space no longer contain any information and the stiff phase-space model can be replaced by a model with only spatial dimensions and without $\epsilon$-dependent time-step constraints. For example, if the steady state distribution density $\mathcal{M}(v)$ is symmetric with mean zero under the Lebesgue measure, in the limit $\epsilon \to 0$, the solution \eqref{eq:kineticdiffusive} can be decomposed as $f(x,v,t) = \rho(x,t)\mathcal{M}(v)$ where $\rho(x,t)$ satisfies the diffusion equation~\cite{Lapeyre2003}
\begin{equation}
\label{eq:heat}
\partial_t \rho(x,t) = \nabla_{xx} \rho(x,t).
\end{equation}
If $\epsilon$ is small enough, the issues of high dimensionality and multiple time-scales can be avoided by approximating \eqref{eq:kineticdiffusive} by the limiting macroscopic equation \eqref{eq:heat}. The closer to the limit, the better the approximation is. Given knowledge of $\epsilon$ over the domain $\D_x$, one could consider a domain decomposition approach, selecting a model based on $\epsilon$, see e.g.~\cite{Illner1993, Klar1995}. In this work, we take an alternate route by using a single model and a corresponding asymptotic-preserving scheme, across the domain.

Asymptotic-preserving schemes (AP)~\cite{Jin2000} incorporate knowledge about the limiting macroscopic equation into the discretization of the kinetic equation. This is done in a way that produces unconditionally stable numerical convergence to the limiting equation as $\epsilon \to 0$, while still maintaining consistency for a fixed finite time-scale separation. For an overview of such methods we refer to~\cite{Dimarco2014, Jin2022}. Many deterministic AP schemes, simulating $f(x,v,t)$, have been developed using various approaches such as operator splitting~\cite{Boscarino2013,  Gosse2002, Jin1999, Jin1998, Jin2000}, truncated expansions in $\epsilon$~\cite{Klar1998,Klar1999}, micro-macro decomposition~\cite{Crouseilles2011,Lemou2008}, spatial averaging~\cite{Larsen1989, Larsen1987} and low-rank approximation~\cite{Einkemmer2021a, Einkemmer2021}. Similar schemes have also been developed for the hyperbolic scaling~\cite{Bennoune2008, Dimarco2011, Dimarco2012, Klar1999a}, i.e., multiplying the right-hand side of \eqref{eq:kineticdiffusive} by $\epsilon$.

Due to the problem's high dimensionality, a stochastic approach with Monte Carlo particle simulation is often favorable. In such an approach, individual particle trajectories are simulated, producing samples from the distribution with density $f(x,v,t)$ whose contribution to the QoI is averaged. By tracing individual particle paths, rather than forming a $2d$-dimensional grid, most of the computational effort is spent sampling high-probability regions of the phase-space. Assuming the QoI is well behaved, e.g., it does not depend disproportionately on low-probability regions of $f(x,v,t)$, this approach will allocate computational resources more efficiently. Monte Carlo methods, however, introduce statistical noise to the computed result, as $f(x,v,t)$ is approximated by a finite number of random samples. A number of stochastic AP schemes have been developed in the hyperbolic scaling, see e.g.~\cite{Degond2011,Dimarco2008,Dimarco2010,Pareschi1999,Pareschi2001,Pareschi2005}. For AP Monte Carlo schemes in the diffusive scaling, we refer to~\cite{Crestetto2018,Dimarco2018,Mortier2019}.

In~\cite{Loevbak2020a, Loevbak2021}, the multilevel Monte Carlo method (MLMC)~\cite{Giles2008,Heinrich2001} was applied to the AP scheme in~\cite{Dimarco2018}, in an approach called asymptotic-preserving multilevel Monte Carlo (APMLMC). The key idea was to leverage the ability of AP schemes to perform stable simulations with a collision-rate-independent time step size to first compute a biased, but low cost estimate of the QoI. By using large time steps, many trajectories can be simulated, so the computed estimate will have low variance. Then, the initial estimate is improved upon with a sequence of corrective estimators for the bias, formed by comparing simulations with different time step sizes. Under the condition that the sum of these corrections converges sufficiently fast, fewer samples will be required as the simulation accuracy is increased. For a detailed description of this approach and an analysis of the method, we refer the reader to~\cite{Loevbak2021}.  Aside from the APMLMC approach, prior work on applying MLMC to particle methods for kinetic equations can be found in \cite{Haji-Ali2018, Rosin2014} as well as extensions in, e.g., importance sampling~\cite{BenRached2022} and iterative methods for treating coupled particles~\cite{Belomestny2019, Szpruch2019}.

The works~\cite{Loevbak2020a, Loevbak2021} demonstrated that the APMLMC approach can achieve a significant speedup when a high accuracy is demanded of the QoI. However, the majority of the method's computational cost lies in simulations with a time step size $\Delta t = \mathcal{O}(\epsilon^2)$, meaning further developments are needed to remove the collision-rate dependence of the simulation cost. In this work, we extend the method in~\cite{Loevbak2020a, Loevbak2021} by introducing a new algorithm for correlating simulations with very coarse time steps and time steps $\Delta t = \mathcal{O}(\epsilon^2)$. A preliminary version of this algorithm was published as~\cite{Loevbak2020}, where numerical experiments showed a significant speedup at the cost of an extra bias. In this paper, we expand on the work of~\cite{Loevbak2020}. We modify the algorithm to reduce the remaining bias, removing it completely for Gaussian steady-state density $\mathcal{M}(v)$. We also perform a detailed numerical analysis of the new algorithm.  We also refer to~\cite{Mortier2020}, where similar ideas were used to develop an alternate APMLMC scheme based on the AP scheme in~\cite{Mortier2019}.

The remainder of this paper is structured as follows. In Section~\ref{sec:kinetic_equations}, we demonstrate the numerical issues that arise in forward Monte Carlo simlulation of kinetic equations near the diffusive limit and present the general idea behind the APMLMC approach. In Section~\ref{sec:improvedcoupling}, we present our improved APMLMC scheme. In Section~\ref{sec:analysis}, we compute expressions for the optimal method parameters and demonstrate why the new scheme outperforms previous work. In Section~\ref{sec:experiments}, we demonstrate both the speedup and correctness achieved by using our new correlation with numerical experiments. Finally, in Section~\ref{sec:conclusions}, we draw some conclusions and discuss possible future work.

\section{Multilevel Monte Carlo Simulation near the Diffusive Limit}
\label{sec:kinetic_equations}
Further in this work, we consider equations of the form \eqref{eq:kineticdiffusive} where the velocity density $\mathcal{M}(v)$ is symmetric with characteristic velocity $\CV$, which we define as the standard-deviation of the corresponding distribution. To simplify later notation, we introduce the probability density $\mathcal{B}(\bar{v})$, with mean 0 and variance 1, and a decomposition
\begin{equation}
\label{eq:velocitydecomp}
\mathcal{M}(v) = \frac{1}{\CV} \mathcal{B}\left( \frac{v}{\CV} \right) = \frac{1}{\CV} \mathcal{B}(\bar{v})
\end{equation}
The two concrete cases we consider for $\mathcal{M}(v)$ are:
\begin{itemize}
\item The Goldstein-Taylor model where $v \in \{-\CV,\CV\}$, i.e., $\mathcal{M}(v) = \frac{1}{2} \left( \delta_{v,-\CV} + \delta_{v,\CV}  \right)$, with $\delta$ the Kronecker delta function.
\item Continuous, normally distributed velocities $v\in\R$, i.e., $\mathcal{M}(v) = \mathcal{N}(v;0,\CV^2)$.
\end{itemize}
One could also consider other velocity distributions, e.g., due to a field driving the particles or a domain boundary. We leave the discussion of non-symmetric distributions and other generalizations for future work.

We consider a quantity of interest $Y(t^*)$, which is an integral over a function $F(x,v)$ of a particle's position $X(t)$ and velocity $V(t)$ at time $t=t^*$, with respect to the measure $f(x,v,t)\,\text{d}x\,\text{d}v$:
\begin{equation}
\label{eq:QoI}
Y(t^*) = \E \left[F\!\left(X(t^*),V(t^*)\right)\right] = \int_{\D_v}\int_{\D_x} F(x,v) f(x,v,t^*)\,\dx\,\dv.
\end{equation}
In applications, the quantity of interest is not always integrated over space. We, however, perform this integration to simplify notation and reduce the cost of simulations.

In Section~\ref{sec:standardMC}, we first demonstrate why straightforward Monte Carlo requires fine time steps to simulate \eqref{eq:kineticdiffusive}. In Section~\ref{sec:ap}, we then present an asymptotic-preserving scheme which allows for simulation with coarse time steps. Section~\ref{sec:mlmc} describes how the multilevel Monte Carlo method can be applied to such a scheme. In Section~\ref{sec:apmlmc}, we describe the APMLMC approach from~\cite{Loevbak2020a, Loevbak2021}, which we improve upon in this work. Finally, in Section~\ref{sec:strategy}, we summarize the level selection strategy used in~\cite{Loevbak2021}.

\subsection{Standard Monte Carlo}
\label{sec:standardMC}

The particle probability density $f(x,v,t)$  at time $t = n \Delta t$ is represented by an ensemble of $P$ particles, with indices $p \in \{1,\dots,P\}$,
\begin{equation}
\label{eq:ensemble}
\left\{\left(X_{p,\Delta t}^n,V_{p,\Delta t}^n\right)\right\}_{p=1}^P.
\end{equation}
We simulate \eqref{eq:kineticdiffusive} with a fixed time-step particle scheme with time step size $\Delta t$. Each particle $p$ has a state $\left( X_{p,\Delta t}^n , V_{p,\Delta t}^n \right)$ in the position-velocity phase space at each time step $n$, with $X_{p,\Delta t}^n \approx X_p(n\Delta t)$ and $V_{p,\Delta t}^n \approx V_p(n\Delta t)$. A classical Monte Carlo estimator $\hat{Y}(t^*)$ for $Y(t^*)$ averages over an ensemble~\eqref{eq:ensemble}
\begin{equation}
\label{eq:MCestimator}
\hat{Y}(t^*) = \frac{1}{P}\sum_{p=1}^P F\!\left(X^N_{p,\Delta t}, V^N_{p,\Delta t}\right), \quad t^* = N \Delta t.
\end{equation}

We use operator splitting to produce a first order scheme in the time step size $\Delta t$~\cite{Pareschi2005a}. For \eqref{eq:kineticdiffusive}, this means alternating between two actions in each time step $n$:
\begin{enumerate}
\item \textbf{Transport step.} Update the particle position given its velocity $V_{p,\Delta t}^n$, i.e.,
\begin{equation}
\label{eq:transport_step}
X^{n+1}_{p,\Delta t} = X^n_{p,\Delta t} + \Delta t V_{p,\Delta t}^n.
\end{equation}
\item \textbf{Collision step.} Perform a collision with probability $p_{c,\Delta t} = \Delta t/\epsilon^2$ by sampling $\mathcal{M}(v)$, i.e.,
\begin{equation}
\label{eq:collision_step}
V_{p,\Delta t}^{n+1} = 
\begin{cases}
V_{p,\Delta t}^{n,*} \sim \frac{1}{\epsilon} \mathcal{M}(v),&\text{with collision probability } p_{c,\Delta t} = \Delta t/\epsilon^2,\\
V_{p,\Delta t}^n,&\text{otherwise,} 
\end{cases}
\end{equation}
\end{enumerate}
where $\sim$ denotes sampling a distribution, which we often represent by its density.

The scheme \eqref{eq:transport_step}--\eqref{eq:collision_step} has a time-step constraint $\Delta t = \mathcal{O}(\epsilon^2)$ when approaching the limit $\epsilon \to 0$. This constraint results in high simulation costs, despite the well defined limit \eqref{eq:heat}~\cite{Dimarco2018}. The variance of the estimator \eqref{eq:MCestimator} decreases with the number of sampled trajectories as $P^{-1}$. This means that a high individual trajectory cost will constrain the certainty of the estimate $\hat{Y}(t^*)$, given a fixed computational budget as taking time steps $\Delta t = \mathcal{O}(\epsilon^2)$ will limit the number of simulated trajectories.

\subsection{Asymptotic-preserving Monte Carlo scheme}
\label{sec:ap}
To remove this time-step constraint, we replace \eqref{eq:transport_step}--\eqref{eq:collision_step} with an asymptotic-preserving scheme, when simulating the ensemble \eqref{eq:ensemble}. The AP scheme from~\cite{Dimarco2018} substitutes \eqref{eq:kineticdiffusive} with
\begin{equation}
\label{eq:GTmod}
\partial_t f + \frac{\epsilon v}{\epsilon^2+\Delta t} \cdot \nabla_x f = \frac{\CV^2 \Delta t}{\epsilon^2+\Delta t} \cdot \nabla_{xx} f + \frac{1}{\epsilon^2+\Delta t} \left(\mathcal{M}(v)\rho - f \right),
\end{equation}
where the space, velocity and time dependency of $f(x,v,t)$ and $\rho(x,t)$ are omitted, for conciseness. Note that the coefficients of \eqref{eq:GTmod} now explicitly contain the simulation time step $\Delta t$. Using this scheme, $\Delta t$ can be chosen independently of $\epsilon$.

Equation \eqref{eq:GTmod} has the following properties~\cite{Dimarco2018}:
\begin{enumerate}
\item For $\frac{\Delta t}{\epsilon^2} \to \infty$, the solution of \eqref{eq:GTmod} converges to a form $f(x,v,t) = \rho(x,t)\mathcal{M}(v)$, with $\rho(x,t)$ satisfying the diffusion equation \eqref{eq:heat} with a bias $\mathcal{O}\big(\frac{\epsilon^2}{\Delta t}\big)$.
\item For $\frac{\Delta t}{\epsilon^2} \to 0$, the coefficients of \eqref{eq:GTmod}  converge to those of the original kinetic equation \eqref{eq:kineticdiffusive} with a bias $\mathcal{O}\big(\frac{\Delta t}{\epsilon^2}\big)$.
\end{enumerate}
When fixing $\Delta t$, the first property leads to the modified kinetic equation having the same asymptotic limit in $\epsilon$ as the original kinetic equation. When fixing $\epsilon$ and varying $\Delta t$ these properties gives an intuitive interpretation to \eqref{eq:GTmod}: We can interpret the modified equation as a combination of the diffusion equation \eqref{eq:heat} and the original kinetic equation \eqref{eq:kineticdiffusive}, with each contribution weighted in function of $\Delta t$. For large time steps, diffusion dominates over kinetic effects, whereas, for small time steps, kinetic effects are dominant. At the particle level, diffusion corresponds to a Brownian motion, which has no constraints on the time step size, hence the stability of \eqref{eq:GTmod}.

Particle trajectories corresponding with \eqref{eq:GTmod} are generated as follows:
\begin{enumerate}
\item \textbf{Transport-diffusion step.} Update the particle given its velocity $V^n_{p,\Delta t}$ and a Brownian increment with diffusion coefficient $\displaystyle D_{\Delta t}=\frac{\CV^2\Delta t}{\epsilon^2 + \Delta t}$, i.e.,
\begin{equation}
\label{eq:transport}
X^{n+1}_{p,\Delta t} = X^n_{p,\Delta t} +\Delta t V^n_{p,\Delta t}  + \sqrt{2 \Delta t}\sqrt{D_{\Delta t}}\xi^n_p,
\end{equation}
in which we generate $\xi_p^n \sim \mathcal{N}(0,1)$.
\item \textbf{Collision step.} Perform a collision with probability $\displaystyle p_{c,\Delta t} = \frac{\Delta t}{\epsilon^2+\Delta t}$ by sampling $\mathcal{M}_{\Delta t}(v)$ with characteristic velocity $\displaystyle\CV_{\Delta t}=\frac{\epsilon}{\epsilon^2+\Delta t}\CV$, i.e.,
\begin{equation}
\label{eq:collision}
V_{p,\Delta t}^{n+1} = 
\begin{cases}
V_{p,\Delta t}^{n,*} \sim \mathcal{M}_{\Delta t}(v),&\text{with probability } p_{c,\Delta t} = \dfrac{\Delta t}{\epsilon^2+\Delta t},\\
V_{p,\Delta t}^n,&\text{otherwise.} 
\end{cases}
\end{equation}

\end{enumerate}
The scheme \eqref{eq:transport}--\eqref{eq:collision} stably generates trajectories for large $\Delta t$. However, these trajectories are biased with $\mathcal{O}(\Delta t)$ compared to those generated by \eqref{eq:transport_step}--\eqref{eq:collision_step}.

\subsection{Multilevel Monte Carlo}
\label{sec:mlmc}

We combine simulations with different time step sizes through multilevel Monte Carlo (MLMC)~\cite{Giles2008,Heinrich2001}. Multilevel Monte Carlo combines many inaccurate trajectories (low variance) with lower numbers of increasingly accurate trajectories (low bias). A coarse Monte Carlo estimator at level 0 $\hat{Y}_0(t^*)$, with a time step $\Delta t_0$, is given by
\begin{equation}
\label{eq:MC0estimator}
\hat{Y}_0(t^*) = \frac{1}{P_0}\sum_{p=1}^{P_0} F\!\left(X^{N_0}_{p,\Delta t_0},V^{N_0}_{p,\Delta t_0}\right), \quad t^* = N_0 \Delta t_0.
\end{equation}
Estimator \eqref{eq:MC0estimator} is based on a large number of trajectories $P_0$, but has a low computational cost as few time steps $N_0$ are required to reach time $t^*$.

Next, $\hat{Y}_0(t^*)$ is refined upon by a sequence of $L$ difference estimators at levels $\ell=1,\ldots,L$. Each difference estimator uses an ensemble of $P_\ell$ particle pairs
\begin{equation}
\label{eq:MClestimator}
\hat{Y}_\ell(t^*) = \frac{1}{P_\ell}\sum_{p=1}^{P_\ell}\left(F\!\left(X^{N_\ell}_{p,\Delta t_\ell},V^{N_\ell}_{p,\Delta t_\ell}\right)-F\!\left(X^{N_{\ell-1}}_{p,\Delta t_{\ell-1}},V^{N_{\ell-1}}_{p,\Delta t_{\ell-1}}\right)\right)\!,\quad t^* = N_\ell \Delta t_\ell.
\end{equation}
Each correlated particle pair consists of a particle with a fine time step $\Delta t_\ell$ and a particle with a coarse time step $\Delta t_{\ell-1} = M \Delta t_\ell$, with $M$ a positive integer. The particles in each pair undergo correlated simulations, which intuitively can be understood as making both particles follow the same qualitative trajectory for two different simulation accuracies. One can interpret the difference estimator \eqref{eq:MClestimator} as using the fine simulation to estimate the bias in the coarse simulation.

Given a sequence of levels $\ell \in \{0,\dots,L\}$, with decreasing time step sizes, and the corresponding estimators given by \eqref{eq:MC0estimator}--\eqref{eq:MClestimator}, the multilevel Monte Carlo estimator for the quantity of interest $Y(t^*)$ is computed by the sum
\begin{equation}
\label{eq:telescopic}
\hat{Y}(t^*) = \sum_{\ell=0}^{L} \hat{Y}_\ell(t^*).
\end{equation}
Using linearity of the expectation, it is clear that the expected value of the estimator \eqref{eq:telescopic} is the same as that of \eqref{eq:MCestimator} with the finest time step $\Delta t = \Delta t_L$. If the number of simulated (pairs of) trajectories $P_\ell$ decreases sufficiently fast as $\ell$ increases, it can be shown that the multilevel Monte Carlo estimator requires a lower computational cost than a classical Monte Carlo estimator to achieve the same mean square error~\cite{Giles2015}.

\subsection{Term-by-term correlation approach}
\label{sec:apmlmc}

The differences in \eqref{eq:MClestimator} will only have low variance if the simulated paths up to $X^{N_\ell}_{p,\Delta t_\ell}$ and $X^{N_{\ell-1}}_{\Delta t_{p,\ell-1}}$, with time steps related by $\Delta t_{\ell-1}=M\Delta t_{\ell}$ are correlated. Before introducing our new correlation in Section~\ref{sec:improvedcoupling}, we first summarize the correlation introduced in~\cite{Loevbak2020a}. To discuss correlated trajectory pairs, we define a sub-step index $m \in \{0,\dots, M-1\}$, i.e., $X^{n,m}_{p,\Delta t_\ell}\equiv X^{nM+m}_{p,\Delta t_\ell}$. To span a time interval of size $\Delta t_{\ell-1}$, the coarse simulation requires a single time step of size $\Delta t_{\ell-1}$ while the fine simulation requires $M$ time steps of size $\Delta t_{\ell}$: 
\begin{equation}
\label{eq:coupled_simulation}
\begin{dcases}
X^{n+1}_{p,\Delta t_{\ell-1}} = X^{n}_{p,\Delta t_{\ell-1}} + \Delta t_{\ell-1}V_{p,\Delta t_{\ell-1}}^n + \sqrt{2 \Delta t_{\ell-1}} \sqrt{D_{\Delta t_{\ell-1}}} \xi^{n}_{p,\ell-1} \\
X^{n+1,0}_{p,\Delta t_{\ell}} = X^{n,0}_{p,\Delta t_{\ell}} + \sum_{m=0}^{M{-}1} \left( \Delta t_\ell V_{p,\Delta t_\ell}^{n,m} + \sqrt{2 \Delta t_{\ell}} \sqrt{D_{\Delta t_\ell}} \xi^{n,m}_{p,\ell}\right)
\end{dcases},
\end{equation}
with $\xi^{n}_{p,\ell-1}, \xi^{n,m}_{p,\ell} \sim \mathcal{N}(0,1)$, $V_{p,\Delta t_{\ell-1}}^n \sim \mathcal{M}_{\Delta t_{\ell-1}}(v)$ and $V_{p,\Delta t_{\ell}}^{n,m} \sim \mathcal{M}_{\Delta t_\ell}(v)$.

To correlate a pair of simulations with time step sizes $\Delta t_{\ell}$ and $\Delta t_{\ell-1}$, we first perform $M$ independent fine time steps of size $\Delta t_\ell$. Then, the resulting $M$ sets of random numbers are combined into a single set of random numbers for the corresponding coarse time step of size $\Delta t_{\ell-1}$. If these coarse random numbers are distributed as if they had been generated independently, then the coarse simulation statistics will be preserved, while also correlating both paths.

There are two sources of stochastic behavior in scheme \eqref{eq:transport}--\eqref{eq:collision}. On the one hand, a new Brownian increment $\xi_p^n$ is generated at each time step. On the other hand, each particle has a non-zero collision probability in each time step. Collisions cause a new particle velocity $V_{p, \Delta t}^{n+1}$ to be sampled for the next time step. The term-by-term correlation presented in~\cite{Loevbak2020a, Loevbak2021} correlates these sources as two separate phenomena. We briefly summarize this approach here, referring to~\cite{Loevbak2021} for more details.

\subsubsection*{Generating $\xi^{n}_{p,\ell-1}$} Each fine sub-step $m \in \{0,\dots, M{-}1\}$ generates a Brownian increment $\xi^{n,m}_{p,\ell} \sim \mathcal{N}(0,1)$. These are summed and rescaled to have unit variance
\begin{equation}
\label{eq:xicorr}
\xi^{n}_{p,\ell-1} = \frac{1}{\sqrt{M}}\sum_{m=0}^{M{-}1} \xi^{n,m}_{p,\ell} \sim \mathcal{N}(0,1),
\end{equation}
 similar to the approach proposed in~\cite{Giles2008}.

\subsubsection*{Generating $V^{n+1}_{p,{\ell-1}}$} At the end of each fine sub-step $m \in \{0,\dots, M{-}1\}$, a collision occurs with probability $p_{c,\Delta t_\ell}$. This process is simulated by sampling a $\PT^{n,m}_{p,\ell} \sim \mathcal{U}([0,1])$, with $\mathcal{U}([0,1])$ the uniform distribution on $[0,1]$. If
\begin{equation}
\label{eq:collisioncondition}
\PT_{p,\ell}^{n,m} \geq p_{nc,\Delta t_\ell} =1-p_{c,\Delta t_\ell} = \frac{\epsilon^2}{\epsilon^2+\Delta t_\ell}
\end{equation}
then a collision takes place at the end of fine sub-step $m$. There is at least one fine-simulation collision if $\PT^{n,\text{max}}_{p,\ell} = \max_m \PT_{p,\ell}^{n,m} \geq p_{nc,\Delta t_\ell}$. When
\begin{equation}
\label{eq:collisioncorr}
\PT_{p,\ell-1}^n \geq \pncr, \quad \PT_{p,\ell-1}^n = {\left(\PT^{n,\text{max}}_{p,\ell}\right)}^M \sim \mathcal{U}([0,1]),
\end{equation}
a collision also occurs in the simulation at level $\ell-1$. If there is a collision at level $\ell-1$, then the last value $\VD^{n,m,*}_{p,\ell}$ from the fine simulation, with $\VD^{n,m,*}_{p,\ell}\!\sim\mathcal{B}(\bar{v})$ as defined in \eqref{eq:velocitydecomp}, is used to generate a velocity for the next step at level $\ell-1$, i.e.,
\begin{equation}
\label{eq:velocitycorr}
V^{n+1}_{p,{\ell-1}} = \CV_{\Delta t_{\ell-1}} \VD^{n+1}_{p,{\ell-1}}, \quad \VD^{n+1}_{p,{\ell-1}} = \VD^{n,*}_{p,{\ell-1}} = \VD^{n,M-1,*}_{p,\ell}.
\end{equation}
If there is no collision in the coarse simulation, then its velocity remains unchanged.

\subsection{Level strategy}
\label{sec:strategy}
In~\cite{Loevbak2021}, the following level sequence was suggested, after observing that difference estimators with $\Delta t_\ell > \epsilon^2$ will not aid in variance reduction:
\begin{enumerate}
\item At level 0, generate an initial estimate of $\hat{Y}(t^*)$ by simulating with $\Delta t_0 = t^*$.
\item At level 1, perform correlated simulations to $t^*$ using $\Delta t_0 = t^*$ and $\Delta t_1 = \epsilon^2$.
\item Continue with a geometric sequence of levels $\Delta t_l = \epsilon^2M^{1-l}$ for $l>1$ until the requested root mean square error bound $E$ is achieved.
\end{enumerate}
It was then shown in~\cite{Loevbak2021} that an MLMC scheme based on the term-by-term approach achieves an asymptotic speed-up $\mathcal{O}\left(\mse^{-2}\log^2(\mse)\right)$ in the root mean square error bound $\mse$, when compared to a single level simulation for the same $\mse$. This asymptotic result, however, does not say how well the presented approach bridges the gap between very coarse (diffusive) simulations and calculations with time steps in the order of magnitude of the collision rate, capturing the collision dynamics. In fact, experimental results showed only a slight speed up, compared to leaving out level 0 altogether~\cite{Loevbak2020a}. In the current paper we therefore introduce an improved correlation at level 1.

%

\section{An improved APMLMC scheme with combined correlation}
\label{sec:improvedcoupling}

\newlength{\figurewidth}
\setlength{\figurewidth}{0.24\textwidth}
\newlength{\figureheight}
\setlength{\figureheight}{\figurewidth}

\begin{figure}
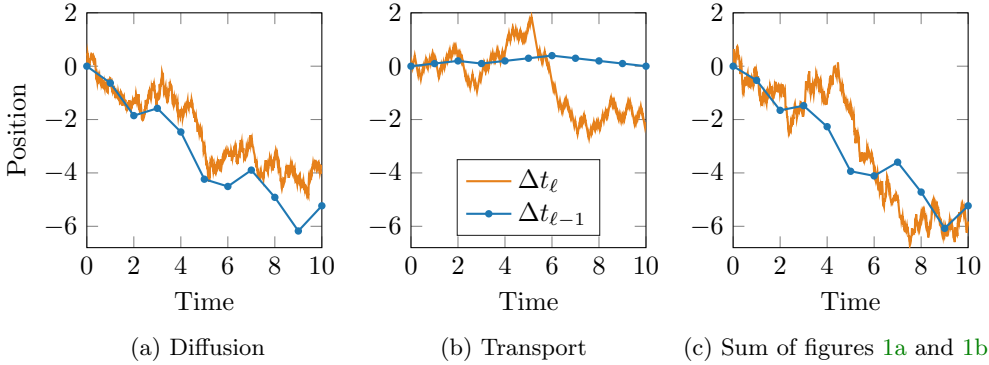

\centering
\begin{subfigure}{0.32\textwidth}
\centering
\hspace{-0.7cm}\input{explainingLevel1/diffusionVisualization}
\caption{Diffusion\label{fig:largeMdiffusion}}
\end{subfigure}
\hspace{-0.2cm}
\begin{subfigure}{0.32\textwidth}
\centering
\input{explainingLevel1/transportVisualization}
\caption{Transport\label{fig:largeMtransport}}
\end{subfigure}
\begin{subfigure}{0.32\textwidth}
\centering
\input{explainingLevel1/combinedVisualization}
\caption{Sum of figures~\ref{fig:largeMdiffusion} and~\ref{fig:largeMtransport}\label{fig:largeMcombined}}
\end{subfigure}
\caption{A pair of trajectories with two discrete velocities and $\CV=1$, correlated by the approach described in Section~\ref{sec:apmlmc}. We set $\epsilon = 0.1$, $\Delta t_\ell = 0.01$ and $\Delta t_{\ell-1} = 1.$\label{fig:largeM}}
\end{figure}

%
\newcommand{\tabletimes}{\!\times\!}
While the term-by-term correlation described in Section~\ref{sec:apmlmc} works well for small refinement factors $M$, it does not produce well correlated trajectory pairs when $M$ is large~\cite{Loevbak2020}, especially so when $\Delta t_\ell \approx \epsilon^2$, as is the case in level 1. The reason for this becomes clear when we consider two trajectories from such a pair, as shown in Figure~~\ref{fig:largeM}, where we separately show the diffusive increments $\sqrt{2\Delta t_{\ell-1}}\sqrt{D_{\Delta t_{\ell-1}}}\xi^{n}_{p,\ell-1}$ and $\sqrt{2\Delta t_{\ell}}\sqrt{D_{\Delta t_\ell}}\xi^{n,m}_{p,\ell}$ (Figure~\ref{fig:largeMdiffusion}), transport increments $\Delta t_{\ell-1} V_{p,\Delta t_{\ell-1}}^{n}$ and $\Delta t_\ell V_{p,\Delta t_\ell}^{n,m}$ (Figure~\ref{fig:largeMtransport}), and the sum of both diffusive and transport increments (Figure~\ref{fig:largeMcombined}). Both trajectories in Figure~\ref{fig:largeMdiffusion} are similar in shape, but have diffusion coefficients $D_{\Delta t}$ that differ by almost a factor two. More noticeably the two trajectories in Figure~\ref{fig:largeMtransport} bare almost no similarity. The reason these trajectories differ are \textbf{(i)} their differing characteristic velocities $\CV_{\Delta t}$; and \textbf{(ii)} that only one of the $M$ fine velocities $\VD_{p,\ell}^{n,m,*}$ is used to generate a coarse velocity $\VD_{p,\ell}^{n,*}$. Note that reason (ii) means that a lot of information in the fine process is lost when generating the coarse process.

Let us start again from \eqref{eq:coupled_simulation} with $\Delta t_{\ell-1} \gg \Delta t_\ell$, bearing in mind the motivation for using an asymptotic-preserving scheme, i.e., substituting a large number of (in the limit $\epsilon \to 0$ infinitely many) collisions with a single diffusive step. From this motivation, it is apparent that there should be a correlation between the fine simulation velocities and the Brownian increment in the coarse simulation, which, for a large part, replaces their behavior. In~\cite{Kloeden1992}, it is shown that one can simulate Brownian motion in the weak sense, using approximate Brownian increments $\xi_{p,\ell-1}^n$ satisfying
\begin{equation}
\label{eq:bound_weak_euler}
\left| \E \left[ \xi_{p,\ell-1}^n \right] \right| + \left| \E \left[ \left( \xi_{p,\ell-1}^n \right)^3 \right] \right| + \left| \E \left[ \left( \xi_{p,\ell-1}^n \right)^2 \right] - 1 \right| \leq K \Delta t_{\ell-1},
\end{equation}
for some constant $K$. If \eqref{eq:bound_weak_euler} is satisfied, the same weak convergence is achieved as the classical Euler-Maruyama scheme. As we are interested in computing moments, having weak convergence is sufficient for maintaining good approximations in coarse simulations. We now present a new scheme for generating approximate Brownian increments $\xi_{p,\ell-1}^n$ consisting of two parts:
\begin{enumerate}
\item We generate a coarse value $\xi_{p,\ell-1}^n$, satisfying \eqref{eq:bound_weak_euler}, through a weighted sum of (a) the fine velocities $\VD_{p,\ell}^{n,m,*}$ and (b) the normaly distributed increments $\xi_{p,\ell}^{n,m}$. The details of the resulting algorithm are presented in Section~\ref{sec:improvedcorrelation}.
\item The coarse values $\xi_{p,\ell-1}^n$ generated by the new algorithm will be normally distributed if the steady state velocity probability density $\mathcal{M}(v) \equiv \mathcal{N}(v;0,1)$. In the general case for $\mathcal{M}(v)$, we can however only assume \eqref{eq:bound_weak_euler} holds. As a result, \eqref{eq:telescopic} will no longer be consistent as the distribution of the values $\xi_{p,\ell-1}^n$ used to generate $\hat{Y}_{\ell}$ no longer match those used to generate $\hat{Y}_{\ell-1}$. In Section~\ref{sec:level0}, we present a modification to the generation of the values $\xi_{p,\ell-1}^n$ in the fine simulation at level $\hat{Y}_{\ell-1}$ which restores consistency in \eqref{eq:telescopic}.
\end{enumerate}
%
%

\subsection{Combined correlation}
\label{sec:improvedcorrelation}

The improved correlation generates $\xi_{p,\ell-1}^n$ as a weighted sum of a diffusion contribution $\xi_{p,\ell-1,W}^n$ from the fine simulation diffusion and a transport contribution $\xi_{p,\ell-1,T}^n$ from the fine simulation rescaled velocities $\VD^{n,m}_{p,\ell}$
\begin{equation}
\label{eq:weighted_sum}
\xi_{p,\ell-1}^n = \sqrt{\theta_{\ell}} \, \xi_{p,\ell-1,W}^n + \sqrt{1-\theta_{\ell}} \, \xi_{p,\ell-1,T}^n,
\end{equation}
with $\theta_{\ell} \in [0,1]$. The physical interpretation of the contribution $\xi_{p,\ell-1,T}^n$ comes from considering that the convergence of \eqref{eq:kineticdiffusive} to \eqref{eq:heat}, at the particle level, means substituting an infinite number of kinetic increments with a single Brownian increment. As the fine simulation has more weight on its kinetic behavior, while the coarse simulation has more weight on its diffusive behavior, basing coarse diffusion on transport increments therefore insures that both simulations follow a similar trajectory.

If the diffusive contribution $\xi_{p,\ell-1,W}^n$ and transport contribution $\xi_{p,\ell-1,T}^n$ are both distributed symmetrically with unit variance, the same holds for $\xi_{p,\ell-1}^n$, meaning that we can set $K \equiv 0$ in \eqref{eq:bound_weak_euler}.
The correlation for generating coarse collisions and velocities remains as in \eqref{eq:collisioncorr}--\eqref{eq:velocitycorr}. In Sections~\ref{sec:combined_diffusion}--\ref{sec:theta}, we discuss the contributions $\xi_{p,\ell-1,W}^n$, $\xi_{p,\ell-1,T}^n$ and the weight $\theta_\ell$, one at a time.

\subsubsection{Diffusive contribution $\xi_{p,\ell-1,W}^n$}\label{sec:combined_diffusion} The coupling of Brownian increments is identical to \eqref{eq:xicorr} from the term-by-term correlation
\begin{equation}
\label{eq:diff_correlation}
\xi_{p,\ell-1,W}^n = \frac{1}{\sqrt{M}} \sum_{m=0}^{M{-}1} \xi_{p,\ell}^{n,m}.
\end{equation}

\subsubsection{Transport contribution $\xi_{p,\ell-1,T}^n$}
\label{sec:combined_transport}
We can generate a value $\xi_{p,\ell-1,T}^n$ with expected value zero and unit variance from the $\VD_{p,\ell}^{n,m}$ as
\begin{equation}
\label{eq:transp_correlation}
\xi_{p,\ell-1,T}^n = \left( \V \left[ \sum_{m=0}^{M{-}1} \VD_{p,\ell}^{n,m} \right] \right)^{-\frac{1}{2}} \sum_{m=0}^{M{-}1} \VD_{p,\ell}^{n,m}.
\end{equation}
We compute the variance of the sum of $\VD_{p,\ell}^{n,m}$ as
\begin{align}
\V \left[ \sum_{m=0}^{M{-}1} \VD_{p,\ell}^{n,m} \right] &= \sum_{m=0}^{M{-}1} \V \left[ \VD_{p,\ell}^{n,m} \right] + 2 \sum_{m=0}^{M-2}\sum_{m^\prime=m+1}^{M{-}1} \Cov\left( \VD_{p,\ell}^{n,m},  \VD_{p,\ell}^{n,m^\prime} \right) \\
&= M + 2 \sum_{\dm=1}^{M{-}1}(M-\dm) \left(\pncf\right)^\dm \\
&= M + 2\frac{\pncf\left((\left(\pncf\right)^M-M(\pncf-1)-1\right)}{(\pncf-1)^2} \\
&= M + 2\frac{\pncf\left(\left(\pncf\right)^M+M\pcf-1\right)}{\left(\pcf\right)^2},\label{eq:variance_beta_sum}
\end{align}
where $\pcf$ is the collision probability as defined in \eqref{eq:collision} and $\pncf=1-\pcf$. 

\begin{remark}
\label{rem:extra_correlation}
As the $\VD^{n,m}_{p,\ell}$ are now used to generate values for both $\xi^n_{p,\ell-1}$ and $\VD^n_{p,\ell-1}$, care must be taken to avoid introducing correlations between the generated $\xi^n_{p,\ell-1}$ and $\VD^n_{p,\ell-1}$. Such correlations can be introduced in two ways:
\begin{enumerate}
\item If there is a collision in both processes during time step $n-1$ then $\VD^{n}_{p,\ell-1} = \VD^{n,0}_{p,\ell}$, by \eqref{eq:velocitycorr}. As a consequence $\VD^{n,0}_{p,\ell}$ and any subsequent $\VD^{n,m}_{p,\ell}$ which fall before the first collision in the process with time step $\Delta t_\ell$ cannot be used in the sum \eqref{eq:transp_correlation}, without introducing a correlation between $\xi^n_{p,\ell-1}$ and $\VD^n_{p,\ell-1}$.
\item If there is no collision in fine time step $n,M-1$ then $\VD^{n+1}_{p,\ell-1} = \VD^{n,M-1}_{p,\ell}$, by \eqref{eq:velocitycorr}. As a consequence any $\VD^{n,m}_{p,\ell}$ which fall after the last collision in the process with time step $\Delta t_\ell$ cannot be used in the sum \eqref{eq:transp_correlation}, without introducing a correlation between $\xi^{n}_{p,\ell-1}$ and $\VD^{n+1}_{p,\ell-1}$.
\end{enumerate}
\end{remark}
\begin{remark}
\label{rem:replaced_velocities}
A preliminary version of this scheme was presented, without analyisis, in~\cite{Loevbak2020}. The scheme presented there suffered from an observable bias due to the dependencies described in Remark~\ref{rem:extra_correlation}. Here, we avoid introducing such dependencies by sampling two values from $\mathcal{M}(v)$: One to substitute the velocities in \eqref{eq:transp_correlation} before the first fine simulation collision, the other to substitute the velocities after the final fine simulation collision during the coarse time interval $n$. Note that no velocities need to be substituted after a final collision taking place in fine sub-step $M-1$.
\end{remark}

\subsubsection{Contribution weight $\theta_{\ell}$}
\label{sec:theta}
What remains is to choose the weight between the two contributions in \eqref{eq:weighted_sum} so as to maximize correlation between the coarse and fine paths.
In Section~\ref{sec:analysis} we analytically derive the optimal value for $\theta_\ell$, for now we simply give a rule-of-thumb: If $M$ is small, as in levels $\ell>1$, then both coarse and fine processes simulate similar models. We can then expect $\theta_\ell$ values close to 1 to work well. If $M$ is large, as in level $\ell=1$, then fine simulation diffusion will have much less effect than diffusion in the coarse simulation, whereas the inverse is true for the transport effects. In this case $\theta_\ell$ values closer to 0 work better as they increase the contribution of fine transport effects to the coarse diffusion.

\subsubsection{Full combined correlation algorithm}
A full simulation of a pair of particles correlated with the combined correlation approach is given in Algorithm~\ref{alg:combined}.

\begin{algorithm}
\caption{Simulating particle pairs with combined correlation.\label{alg:combined}}
\begin{algorithmic}
\FOR{Each time step $n$}
\FOR{$m = 0 \dots M-1$}
\STATE Simulate \eqref{eq:transport}--\eqref{eq:collision} with $\Delta t_\ell$, saving the $\xi^{n,m}_{p,\ell}$, $\PT^{n,m}_{p,\ell}$ and $\VD^{n,m}_{p,\ell}$.
\ENDFOR
\STATE Generate $\xi^{n}_{p,\ell-1}$ from the $\xi^{n,m}_{p,\ell}$ and the $\VD^{n,m}_{p,\ell}$ using \eqref{eq:weighted_sum} and Remark~\ref{rem:replaced_velocities}.
\STATE Generate $\PT_{p,\ell-1}^n$ and $\VD_{p,{\ell-1}}^{n,*}$ from the $\PT_{p,\ell}^{n,m}$ and $\VD_{p,{\ell}}^{n,m,*}$ using \eqref{eq:collisioncorr}--\eqref{eq:velocitycorr}.
\IF{$\PT_{p,\ell-1}^n \geq p_{nc,\Delta t_{\ell-1}}$}
\STATE Set $\VD_{p,{\ell-1}}^{n+1} = \VD_{p,{\ell-1}}^{n,*}$.
\ELSE
\STATE Set $\VD_{p,{\ell-1}}^{n+1} = \VD_{p,{\ell-1}}^{n}$.
\ENDIF
\ENDFOR
\end{algorithmic}
\end{algorithm}

\subsection{Maintaining consistency in the telescopic sum}
\label{sec:level0}

The weak diffusive process produced by \eqref{eq:weighted_sum} does not have the same distribution as the Brownian process used in the independent simulation of \eqref{eq:transport}--\eqref{eq:collision} at level 0, unless $\mathcal{M}(v) \equiv \mathcal{N}(v;0,1)$. A similar issue arises in the alternate scheme in~\cite{Mortier2020}. For large $M$, one can assume that the law of large numbers will result in a suitably small bias. It is however hard to draw any conclusions about the size the error for general quantities of interest due to this inconsistency. We therefore modify the independent level 0 simulation by generating independent coarse simulation Brownian increments $\xi_{p,0}^n$ which are distributed as though they were generated using Algorithm~\ref{alg:combined}, without performing $\mathcal{O}(M)$ computation. To this end, we approximate the right-hand sum in \eqref{eq:transp_correlation}.

First, we consider that \eqref{eq:transp_correlation} does not contain $M$ independent values, as subsequent values of $\VD_{p,\ell}^{n,m}$ will be identical in the absence of a collision. This means that we can rewrite the right-hand sum in \eqref{eq:transp_correlation} as a sum of $I_n \leq M$ independent velocities $\VD_{p,\ell}^{n,i}, i=1,\dots,I_n$, which are maintained for $\lambda_i>0$ subsequent time steps. Here, $\lambda_i \sim \mathcal{G} \left( \pncf \right)$ for $i=1,\dots,I_n-1$, with $\mathcal{G}(\pncf)$ the geometric distribution corresponding with the number of time steps of size $\Delta t_\ell$ required for a collision to occur, and $\lambda_{I_n}$ is chosen so that $\sum_{i=1}^{I_n} \lambda_i \equiv M$. Note that $I_n$ is not fixed, but takes a different value at each time step $n$. We now rewrite the right-hand sum of \eqref{eq:transp_correlation} as
\begin{equation}
\label{eq:rewrite_transport_sum1}
\sum_{m=0}^{M{-}1} \VD_{p,\ell}^{n,m} = \sum_{i=1}^{I_n} \lambda_i \VD_{p,\ell}^{n,i}.
\end{equation}

Next, we group the right-hand side terms of \eqref{eq:rewrite_transport_sum1} by their value of $\lambda_i$, i.e.,
\begin{equation}
\label{eq:rewrite_transport_sum2}
\sum_{i=1}^{I_n} \lambda_i \VD_{p,\ell}^{n,i} = \sum_{\lambda^*=1}^M \sum_{i|\lambda_i = \lambda^*} \lambda_i \VD_{p,\ell}^{n,i} = \sum_{\lambda^*=1}^M \lambda^* \sum_{i|\lambda_i = \lambda^*} \VD_{p,\ell}^{n,i} \approx \sum_{\lambda^*=1}^\Lambda \lambda^* \sum_{i|\lambda_i = \lambda^*} \VD_{p,\ell}^{n,i},
\end{equation}
where we make the conceptual distinction between the sampled random numbers $\lambda_i$ and the index $\lambda^*$.
In the right-hand side approximation of \eqref{eq:rewrite_transport_sum2}, we assume that long sequences of time steps without intermediate collisions are sufficiently rare that they can be assumed not to take place, i.e., given a suitable value $\Lambda < M$, there is no $i$ so that $\lambda_i = \lambda^*$ with $\lambda^* > \Lambda $. This assumption is reasonable at level 1, where $\Delta t_0 \gg \epsilon^2$.

We now define the probability distribution $\mathcal{S}\left( \phi \right)$, representing the sum of $\phi$ independent velocities $\VD_{p,\ell}$ sampled from the velocity distribution with density $\mathcal{M}\left( v \right)$. With this distribution, we can rewrite the approximate right-hand side of \eqref{eq:rewrite_transport_sum2} as
\begin{equation}
\label{eq:rewrite_transport_sum3}
\sum_{\lambda^*=1}^\Lambda \lambda^* \sum_{i|\lambda_i = \lambda^*} \VD_{p,\ell}^{n,i} = \sum_{\lambda^*=1}^\Lambda \lambda^* \sigma_{\lambda^*}, \quad \sigma_{\lambda^*} \sim \mathcal{S} \left( \phi_{\lambda^*} \right),
\end{equation}
with $\phi_{\lambda^*}$ representing the number values $i$ for which $\lambda_i=\lambda^*$. The right-hand side of \eqref{eq:rewrite_transport_sum3} provides us with an approximation of \eqref{eq:transp_correlation} with a computational cost which depends linearly on $\Lambda$ and no longer depends on $M$ assuming the following:
\begin{enumerate}
\item We are able to generate a set of values $\phi_{\lambda^*}$ which are distributed as if we had simulated $M$ potential collisions and counted the number of runs of length $\lambda^*$, under the constraint that no runs longer than $\Lambda$ can occur. Generating these values must have a computational cost independent of $M$.
\item We are able to sample $\mathcal{S}\left( \phi_{\lambda^*} \right)$ with a computational cost independent of $\phi_{\lambda^*}$.
\end{enumerate}
We address how to perform the two points above in the following two subsections. The quality of the approximation also depends directly on $\Lambda$. For $\Lambda = M$, \eqref{eq:rewrite_transport_sum3} is exact, but as $\Lambda$ is decreased then longer sequences without collisions are ignored. In practice a suitable value for $\Lambda$ must be selected, so that the bias resulting from the truncation is negligible compared to the bias of the full multilevel simulation.

\subsubsection{Sampling \texorpdfstring{\boldmath$\phi_{\lambda^*}$}{\textphi\textsubscript{\textlambda\textsuperscript{*}}} from the run-length distribution}
\label{sec:runlengthsampling}
In this section we address how to sample the distribution of the number of runs of $\lambda^*$ time steps with the same velocity, due to no collision taking place. That is, we wish to generate a set of values $\phi_{\lambda^*}$ for $\lambda^* \in \{1,\dots,\Lambda\}$ that are distributed as though we had sequentially simulated $M$ fine time steps with a collision probability $\pncf$ and retained a count of the number of times $\lambda^*$ subsequent time steps maintain the same velocity without intermediate collisions, while enforcing a maximal length $\Lambda$ for such sequences. It is not possible to sample each $\phi_{\lambda^*}$ independently for a time step $n$, given the constraint
\begin{equation}
\label{eq:run_length_constraint}
\sum_{\lambda^*=1}^\Lambda \lambda^* \phi_{\lambda^*} = M.
\end{equation}

In general, it is not straightforward to derive closed-form expressions for the distributions of such runs, however in \cite{Fu1994} a practical approach was developed for sampling various run length statistics of Bernoulli trials making use of Markov chain embedding. To use  this work, we rely on the fact that a run of $\lambda^*$ time steps without intermediate collisions, means that $\lambda^*{-}1$ sequential time steps have been simulated without a collision and that there was a collision in the last of these $\lambda^*$ time steps. This means that we can consider a Bernoulli process on the $M-1$ boundaries between the $M$ fine time steps spanning the coarse time step $\Delta t_{\ell-1}$, where a success is defined as no collision taking place and a failure means that a collision has taken place.

\begin{figure}
\centering
\begin{tikzpicture}[very thick, auto]
\tikzstyle{h} = [circle, draw=blue, fill=blue, minimum width=5pt, fill, inner sep=0pt]
\tikzstyle{m} = [circle, draw=red, minimum width=5pt, inner sep=0pt]
\node[h] (0) at (0,0) {};
\node[m] (1) at (1,0) {};
\node[m] (2) at (2,0) {};
\node[h] (3) at (3,0) {};
\node[m] (4) at (4,0) {};
\node[m] (5) at (5,0) {};
\node[m] (6) at (6,0) {};
\node[h] (7) at (7,0) {};
\node[m] (8) at (8,0) {};
\node[m] (9) at (9,0) {};
\node[m] (10) at (10,0) {};
\node[m] (11) at (11,0) {};
\node[h] (12) at (12,0) {};
\end{tikzpicture}
\caption{A Bernoulli process: Blue dots denote collision, red circles denote no-collision.\label{fig:bernoulli}}
\end{figure}
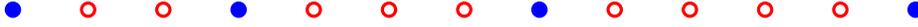

On this Bernoulli process we then define $E_{M-1,\lambda^*{-}1}$ to be the number of success runs of size exactly $\lambda^*{-}1$, i.e., the number of times there are exactly $\lambda^*{-}1$ subsequent time steps without a collision. We also define $G_{M-1,\lambda^*_1}$ to be the number of success runs of size $\lambda^*{-}1$ or larger, i.e., the number of times there are at least $\lambda^*{-}1$ sequential time steps without a collision. In both cases, two success runs are always divided by at least one collision. As an example consider Figure~\ref{fig:bernoulli}, where $M=14$. If we consider $\lambda^*=3$, then $E_{M-1,\lambda^*{-}1} = 1$ and $G_{M-1,\lambda^*{-}1} = 3$. Using the work in~\cite{Fu1994}, we precompute tables for the distributions of $E_{M-1,\lambda^*{-}1}$ and $G_{M-1,\lambda^*{-}1}$ for $\lambda^* = 2, \dots, \Lambda$ for a fixed value of $M$. These tables consist of rows containing the cumulative probabilities of subsequent values of $E_{M-1,\lambda^*{-}1}$ and $G_{M-1,\lambda^*{-}1}$ for a fixed value of $M$, which can then be used for inverse transform sampling in each time step. For the technical details on how these tables are generated, we refer to Appendix~\ref{app:run_length}. While a na\"ive implementation of the algorithm presented in Appendix~\ref{app:run_length} scales as $\mathcal{O}(M^3)$, in practice the time needed for this precomputation proves negligible in comparison to the runtimes of the experiments presented in Section~\ref{sec:experiments}.

Now that we can sample the distributions of runs of no collisions $E_{M-1,\lambda^*{-}1}$ and $G_{M-1,\lambda^*{-}1}$ for a given value of $M$ and $\lambda^*$, we are able to design an algorithm to compute a set of runs of time steps. If we were only interested in computing the number of runs of time steps of a certain length $\lambda^*$, it would be sufficient to simply sample $E_{M-1,\lambda^*{-}1}$ to estimate this result. We, however, want to generate a set of these time-step run-lengths for multiple values of $\lambda^*$, where \eqref{eq:run_length_constraint} needs to hold. This means that we need to sample the conditional probability of $E_{M-1,\lambda^*{-}1}$, given the already sampled numbers. If we generate the disparate $\phi_{\lambda^*}$ in order of decreasing $\lambda^*$ we can consider two cases. If all previously generated $\phi_{\lambda^*}$ were 0, then this conditional probability is equal to the unconditional case and we can use $E_{M-1,\lambda^*{-}1}$ for sampling. If a non-zero $\phi_{\lambda^*}$ has already been generated then we approximate this conditional probability by sampling $G_{M-1,\lambda^*{-}1}$ for the given value of $\lambda^*$ and subtracting the sum of all $\phi_{\lambda^*}$ generated previously. Once all $\phi_{\lambda^*}$ with $\lambda^* > 1$ have been generated, $\phi_1$ is set so that \eqref{eq:run_length_constraint} holds. This process is shown in Algorithm~\ref{alg:run_lengths}.

\begin{algorithm}
\caption{Sampling the distribution of run-lengths of time steps}
\label{alg:run_lengths}
\begin{algorithmic}
\STATE{$\texttt{run\_sum} \gets 0, \quad \texttt{remaining\_steps} \gets M$}
\FOR{$\lambda^* = \Lambda, \dots, 2$}
\STATE{$u \sim \mathcal{U}([0,1])$}
\IF{$\texttt{run\_sum} == 0$}
\STATE{$\phi_{\lambda^*} \gets \texttt{sample\_E\_table}(\lambda^*, u)$}
\ELSE
\STATE{$\phi_{\lambda^*} \gets \texttt{sample\_G\_table}(\lambda^*, u)$}
\ENDIF
\STATE{$\phi_{\lambda^*} \gets \max(\phi_{\lambda^*}-\texttt{run\_sum},0)$}
\STATE{$\texttt{remaining\_steps} \gets \texttt{remaining\_steps} - \lambda^* \phi_{\lambda^*}$}
\WHILE{$\texttt{remaining\_steps} < 0$}
\STATE{$\texttt{remaining\_steps} \gets \texttt{remaining\_steps} + \lambda^*$}
\STATE{$\phi_{\lambda^*} \gets \phi_{\lambda^*} - 1$}
\ENDWHILE
\STATE{$\texttt{run\_sum} \gets \texttt{run\_sum} + \phi_{\lambda^*}$}
\ENDFOR
\STATE{$\phi_1 \gets \texttt{remaining\_steps}$} 
\end{algorithmic}
\end{algorithm}

Algorithm~\ref{alg:run_lengths} does not sample the exact run-length distribution for two reasons:
\begin{enumerate}
\item Run-lengths $\lambda^* > \Lambda$ are not generated in this approach. For sufficiently large $\Lambda$, the occurrence of such run-lengths is rare, but they can strongly impact the variance of the resulting sum of velocities when they occur. This means that a good value for $\Lambda$ is needed to make the trade-off between reducing the computational cost at level 0 and not adding a large bias to the simulation.
\item We do not sample the true conditional probability of run-lengths, subject to the longer run-lengths which have already been sampled. This is because it is possible that the value returned by sampling $G_{M-1,\lambda^*{-}1}$ may be smaller than the sum of $\phi_{\lambda^*}$ values generated up to that point. If this is the case, we set $\phi_{\lambda^*}=0$ for the current $\lambda^*$, which introduces a bias towards longer run-lengths in the sampled distribution. As the expected value of $E_{M-1,\lambda^*{-}1}$ increases with decreasing $\lambda^*$, this is however not a common occurrence. It is also likely that this bias is compensated in part by the fact that the value of $\phi_{\lambda^*}$ is set so that \eqref{eq:run_length_constraint} is satisfied, i.e., the number of intermediate run lengths will on average be slightly too small, while the number of run-lengths of size 1 will be slightly to large to compensate this.
\end{enumerate}

\subsubsection{Sampling sums of velocities}
Once a list of values for $\phi_{\lambda^*}$ has been generated in the case of $\mathcal{M}(v) \equiv \mathcal{N}(v;0,1)$, we need to sample the distribution $\mathcal{S}\left( \phi_{\lambda^*} \right)$. In the general case, there is no closed form expression for sampling this distribution at a cost independent of $\phi_{\lambda^*}$. We solve this issue by generating a table for the two-velocity model. In general, such a table can be generated for any model with a discrete set of velocities, but as the number of velocities grows, so does the size of the table.

\begin{figure}
\centering
\resizebox{0.9\textwidth}{!}{
\begin{tikzpicture}[
                   grow = right,
edge from parent/.style = {draw,-latex},
         label distance = 2mm,
      every node/.style = {minimum width=1.5em, inner sep=2pt},
         level distance = 29mm,
       sibling distance = 12mm,
       baseline]
\node[draw=none,label=90:{$\phi=0$}] {0}
    child {node {-1} 
        child {node {-2}
            child {node {-3}
            	child {node {-4} edge from parent node[above]{$\frac{1}{16}$}}
                child {node {-2} edge from parent node[above]{$\frac{1}{16}$}}
            edge from parent node[above]{$\frac{1}{8}$}}
            child {node {-1}
            	child {node {} edge from parent node[above]{$\frac{3}{16}$}}
            	child {node {0} edge from parent node[above]{$\frac{3}{16}$}}
            edge from parent node[above]{$\frac{1}{8}$}}
                edge from parent node[above]{$\frac{1}{4}$}}
        child {node {}edge from parent node[above]{$\frac{1}{4}$}}
            edge from parent node[above]{$\frac{1}{2}$}}
    child {node[label=90:{$\phi=1$}] {1}
        child {node {0}
            child {node {}edge from parent node[above]{$\frac{2}{8}$}}
            child {node {1}
               	child {node {} edge from parent node[above]{$\frac{3}{16}$}}
            	child {node {2} edge from parent node[above]{$\frac{3}{16}$}}
            edge from parent node[above]{$\frac{2}{8}$}}
                edge from parent node[above]{$\frac{1}{4}$}}
        child {node[label=90:{$\phi=2$}] {2}
            child {node {}edge from parent node[above]{$\frac{1}{8}$}}
            child {node[label=90:{$\phi=3$}] {3}
               	child {node {} edge from parent node[above]{$\frac{1}{16}$}}
            	child {node[label=90:{$\phi=4$}] {4} edge from parent node[above]{$\frac{1}{16}$}}
            edge from parent node[above]{$\frac{1}{8}$}}
                edge from parent node[above]{$\frac{1}{4}$}}
                edge from parent node[above]{$\frac{1}{2}$}};

\end{tikzpicture}
}
\caption{The distribution $\mathcal{S}(\phi)$ for two discrete velocities $\pm 1$.\label{fig:velocity_sum}}
\end{figure}
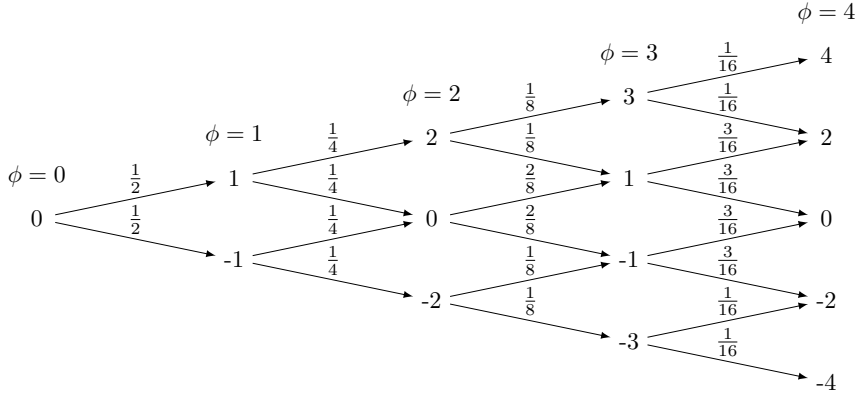

As we sum velocities $\pm 1$, the probability of their sum taking a given value is shown in Figure~\ref{fig:velocity_sum}. It is clear that small velocities, in absolute value, are more probable than larger ones. This means that we can efficiently sample the absolute value of $\sigma_\phi \sim \mathcal{S}(\phi)$ through inverse transform sampling using a table of different potential values of $|\sigma_\phi|$ in increasing order, i.e., sampling a uniformly distributed value $u$, iterating over a the table row corresponding with $\phi$ until a cumulative probability larger than $u$ is encountered, and generating a random sign.

\section{Variance of combined correlation}
\label{sec:analysis}

To demonstrate how the variance of differences of particle positions decreases through the use of the combined correlation, we re-introduce some notation and results from~\cite{Loevbak2021}. The fine and coarse particle position updates after a time step of size $\Delta t_\ell$ or $\Delta t_{\ell-1}$, respectively, are written as
\begin{align}
X_{p,\ell}^{n,m+1} - X_{p,\ell}^{n,m} &= \Delta X_{p,\ell}^{n,m} = \Delta W_{p,\ell}^{n,m} + \Delta T_{p,\ell}^{n,m},\\
X_{p,\ell-1}^{n+1} - X_{p,\ell-1}^n &= \Delta X_{p,\ell-1}^n = \Delta W_{p,\ell-1}^n + \Delta T_{p,\ell-1}^n,
\end{align}
with $\Delta T_{p,\ell}^{n,m} \!=\! \Delta t_\ell V_{p,\Delta t_\ell}^{n,m}$, $\Delta W_{p,\ell}^{n,m} \!=\! \sqrt{2 \Delta t_{\ell}} \sqrt{D_{\Delta t_\ell}} \xi^{n,m}_{p,\ell}$, 
$\Delta T_{p,\ell-1}^n \!=\! \Delta t_{\ell-1} V_{p,\Delta t_{\ell-1}}^n$ and $\Delta W_{p,\ell-1}^n \!=\! \sqrt{2 \Delta t_{\ell-1}} \sqrt{D_{\Delta t_{\ell-1}}} \xi^n_{p,\ell-1}$.
We denote the difference of these updates by $$\displaystyle\Delta_{X,n} = \sum_{m=0}^{M{-}1} \Delta X_{p,\ell}^{n,m} - \Delta X_{p,\ell-1}^n.$$

Using this notation we write out the variance of the difference in position of two simulations after $N$ time steps of size $\Delta t_{\ell-1}$ as
\begin{align}
\label{eq:total_variance}
\V \left[ \sum_{n=0}^{N-1} \Delta_{X,n} \right] = \sum_{n=0}^{N-1} \V \left[  \Delta_{X,n} \right] + 2 \sum_{n=0}^{N-2} \sum_{n^\prime=n+1}^{N-1} \Cov \left( \Delta_{X,n}, \Delta_{X,n^\prime} \right).
\end{align}
As subsequent coarse diffusive increments are independent given the correction in Remark~\ref{rem:replaced_velocities}, we know the covariance term in \eqref{eq:total_variance} is only non-zero due to the dependency between subsequent simulation velocities. As velocities are generated the same way in both combined and term-by-term schemes, we repeat the result computed in~\cite{Loevbak2021}

\begin{align}
&\sum_{n=0}^{N-2} \! \sum_{n^\prime=n+1}^{N-1} \!\!\!\! \Cov(\Delta_{X,n},\Delta_{X,n^\prime}) = \left(\epsilon^2 \CV^2 \left( \pncf^{1-M} + \pncf^{M+1} - 2 \pncf \right) - \Delta t_{\ell-1}^2 \CV_{\Delta t_{\ell-1}}^2 \right) \\
&\qquad\pncf^M \!  \cfrac{\pncf^{MN} \!  {-}  N\pncf^M  \!{+}  N  {-}  1}{\left( 1 {-} \pncf^M \right)^2} + \pncr \!\!  \frac{\pncr^{N} \! {-}  N\pncr \! {+}  N  {-}  1}{\left( 1{-}\pncr \right)^2} \\
&\qquad\qquad\quad\Delta t_{\ell-1}^2\CV_{\Delta t_{\ell-1}}^2 \!\!\! \left( \!\! 1  {-}  \frac{\Delta t_\ell \CV_{\Delta t_\ell}}{\epsilon^2 \CV_{\Delta t_{\ell-1}}}\! \! \left( \! \cfrac{ \pncr{-}\pncf^M }{ 1{-}\pncf^M }\!\!  \left( \!  \pcf^{-1}\! {-} \cfrac{M\pncf^{M}}{ 1{-}\pncf^M}\! \right) \!\! + \! \frac{\epsilon^2}{\Delta t_\ell}\! \right) \!\!\! \right) \label{eq:final_covar}.
\end{align}
The variance term in \eqref{eq:total_variance} is decomposed as
\begin{align}
\V\! \left[  \Delta_{X,n} \right] \!&=\!\! \V\!\! \left[ \sum_{m=0}^{M{-}1}\!\! \Delta X_{p,\ell}^{n,m} \!\!  -\! \! \Delta X_{p,\ell-1}^n \! \right] \!\! = \!\! \V\!\! \left[ \sum_{m=0}^{M{-}1} \!\!\left( \Delta W_{p,\ell}^{n,m} \!\! + \! \Delta T_{p,\ell}^{n,m} \right) \!\! - \!\! \left( \Delta W_{p,\ell-1}^n \!\! + \! \Delta T_{p,\ell-1}^n \right) \! \right]\\
&=\!\! \V\!\! \left[ \sum_{m=0}^{M{-}1} \Delta W_{p,\ell}^{n,m} \right] + \V\!\! \left[ \sum_{m=0}^{M{-}1} \Delta T_{p,\ell}^{n,m} \right] + \V\! \left[ \Delta W_{p,\ell-1}^n \right] + \V\! \left[ \Delta T_{p,\ell-1}^n \right]\label{eq:variance_decomposed}\\
&\quad+ \underbrace{2 \Cov\!  \left( \sum_{m=0}^{M{-}1} \! \Delta W_{p,\ell}^{n,m} , \sum_{m^\prime=0}^{M{-}1} \Delta T_{p,\ell}^{n,m\prime} \right)}_{=0} + \underbrace{2 \Cov\!  \left( \Delta W_{p,\ell-1}^n , \Delta T_{p,\ell-1}^n \right)\vphantom{2 \Cov\!  \left( \sum_{m=0}^{M{-}1} \Delta W_{p,\ell}^{n,m} , \sum_{m^\prime=0}^{M{-}1} \Delta T_{p,\ell}^{n,m\prime} \right)}}_{=0}\\
&\quad- 2 \Cov\!  \left( \sum_{m=0}^{M{-}1} \Delta W_{p,\ell}^{n,m} , \Delta W_{p,\ell-1}^n \right) - \underbrace{2 \Cov\!  \left( \sum_{m=0}^{M{-}1} \Delta W_{p,\ell}^{n,m} , \Delta T_{p,\ell-1}^n \right)}_{=0}\\
&\quad- 2 \Cov\!  \left( \sum_{m=0}^{M{-}1} \Delta T_{p,\ell}^{n,m} , \Delta W_{p,\ell-1}^n \right) - 2 \Cov\!  \left( \sum_{m=0}^{M{-}1} \Delta T_{p,\ell}^{n,m} , \Delta T_{p,\ell-1}^n \right),
\end{align}
where we indicate covariance terms that are zero due to independence. The results
\begin{equation}
\V \left[ \sum_{m=0}^{M{-}1} \Delta W_{p,\ell}^{n,m} \right] = 2 M \Delta t_\ell D_{\Delta t_\ell}, \quad \V \left[ \Delta W_{p,\ell-1}^n \right] = 2 \Delta t_{\ell-1} D_{\Delta t_{\ell-1}} \quad \text{and}
\end{equation}
\begin{align}
&\V \left[ \sum_{m=0}^{M{-}1} \Delta T_{p,\ell}^{n,m} \right] + \V \left[ \Delta T_{p,\ell-1}^n \right]
- 2 \Cov  \left( \sum_{m=0}^{M{-}1} \Delta T_{p,\ell}^{n,m} , \Delta T_{p,\ell-1}^n \right)\\
&\quad\quad = M\Delta t_\ell^2 \CV_{\Delta t_\ell}^2  - \Delta t_{\ell-1}^2  \CV_{\Delta t_{\ell-1}}^2  + 2\epsilon^2 \CV_{\Delta t_\ell}^2  \left( M\Delta t_\ell - (\epsilon^2  +  \Delta t_\ell)\left(1  -  \pncf^M \right)  \right)
\end{align}
can be reused from~\cite{Loevbak2021}, as their derivations remain unchanged. This leaves the two non-zero covariances involving $\Delta W_{p,\ell-1}^n$ to be derived for the new correlation.

The covariance between fine and coarse Brownian updates is straightforward to compute given~\eqref{eq:weighted_sum}--\eqref{eq:diff_correlation}, and the independence and unity variance of the $\xi_{p,\ell}^{n,m} $:

\begin{align}
\Cov\! \left( \sum_{m=0}^{M{-}1}\!\! \Delta W_{p,\ell}^{n,m}\! , \Delta W_{p,\ell-1}^n \!\! \right) \!&= \E\!\left[\! \left( \sum_{m=0}^{M{-}1}\!\!\! \sqrt{2 \Delta t_{\ell}} \sqrt{D_{\Delta t_\ell}} \xi_{p,\ell}^{n,m} \right)\!\! \sqrt{2 \Delta t_{\ell-1}} \sqrt{D_{\Delta t_{\ell-1}}} \xi^n_{p,\ell-1} \right]\\
&=2 \frac{\Delta t_{\ell-1}}{\sqrt{M}} \sqrt{D_{\Delta t_\ell}D_{\Delta t_{\ell-1}}} \E\!\left[ \sum_{m=0}^{M{-}1} \xi_{p,\ell}^{n,m} \xi^n_{p,\ell-1} \right]\\
&=2 \frac{\Delta t_{\ell-1}}{M} \sqrt{\theta_\ell D_{\Delta t_\ell}D_{\Delta t_{\ell-1}}} \E\!\left[  \sum_{m=0}^{M{-}1} \left( \xi_{p,\ell}^{n,m}  \right)^2\right]\\
&=2 \Delta t_{\ell-1}\sqrt{\theta_\ell D_{\Delta t_\ell}D_{\Delta t_{\ell-1}}}.
\label{eq:cov_diffusioncorrelation}
\end{align}

The remaining term in \eqref{eq:variance_decomposed} can be worked out according to \eqref{eq:weighted_sum} and \eqref{eq:transp_correlation} as
\begin{align}
\Cov \!\!  \left(\! \sum_{m=0}^{M{-}1} \!\! \Delta T_{p,\ell}^{n,m}\!\! \!, \Delta W_{p,\ell-1}^n\!\! \right) \!&=\! \E\! \left[ \!\! \left( \! \sum_{m=0}^{M{-}1} \! \Delta t_\ell \CV_{\Delta t_\ell} \bar{V}_{p,\ell}^{n,m}\!\! \right) \!\!\! \left(\!\! \sum_{m^\prime=0}^{M{-}1} \!\! \sqrt{\frac{(1{-}\theta) 2 \Delta t_{\ell-1}\! D_{\Delta t_{\ell-1}}\!}{\V \left[ \sum_{m^\prime=0}^{M{-}1} \VD_{p,\ell}^{n,m^\prime} \right]}} \VD_{p,\ell}^{n,m^\prime}\!\! \right)\!\! \right]\\
= \Delta t_\ell& \CV_{\Delta t_\ell} \sqrt{\frac{(1{-}\theta_\ell) 2 \Delta t_{\ell-1} D_{\Delta t_{\ell-1}}}{\V \left[ \sum_{m^\prime=0}^{M{-}1} \VD_{p,\ell}^{n,m^\prime} \right]}} \E \! \left[ \! \left( \sum_{m=0}^{M{-}1} \bar{V}_{p,\ell}^{n,m} \right) \!\!\! \left( \sum_{m^\prime=0}^{M{-}1}  \VD_{p,\ell}^{n,m^\prime} \right) \! \right]\\
= \Delta t_\ell& \CV_{\Delta t_\ell} \sqrt{\frac{(1{-}\theta_\ell) 2 \Delta t_{\ell-1} D_{\Delta t_{\ell-1}}}{\V \left[ \sum_{m=0}^{M{-}1} \VD_{p,\ell}^{n,m^\prime} \right]}} \sum_{m=0}^{M{-}1} \sum_{m^\prime=0}^{M{-}1}   \E \left[ \VD_{p,\ell}^{n,m} \VD_{p,\ell}^{n,m^\prime} \right]\!.
\label{eq:cov_crosscorrelation}
\end{align} 
The right hand side expected value takes the value 1 if $\VD_{p,\ell}^{n,m}$ and $\VD_{p,\ell}^{n,m^\prime}$ are the same due to the correlation and 0 otherwise. These values are the same if
\begin{enumerate}
\item velocities $m$ and $m^\prime$ have not been replaced as per Remark~\ref{rem:replaced_velocities} and
\item no collision has taken place between the fine time steps $m$ and $m^\prime$.
\end{enumerate}
The veracity of the first point is determined by the index $m$ or $m^\prime$ which is the closest to the beginning or end of the coarse time interval, depending on which case is being considered. On the one hand, we consider whether there is at least one collision between fine sub step 0 and the smaller of $m$ or $m^\prime$. On the other, we consider whether there is at least one collision between the end of the coarse time interval and the larger of $m$ or $m^\prime$. This gives us that the first point is true with probability
\begin{equation}
\label{eq:boundary_probability}
\left( 1 - \pncf^{\min(m,m^\prime)} \right) \left( 1 - \pncf^{M-\max(m,m^\prime)} \right).
\end{equation}
The second point is true if there is no collision in-between fine time steps $m$ and $m^\prime$, which has probability $\pncf^{\left| m-m^\prime \right|}$. Combining both probabilities gives us 
\begin{equation}
\label{eq:single_expectation_crosscorr}
\E \left[ \VD_{p,\ell}^{n,m} \VD_{p,\ell}^{n,m^\prime} \right] = \left( 1 - \pncf^{\min(m,m^\prime)} \right) \left( 1 - \pncf^{M-\max(m,m^\prime)} \right) \pncf^{\left| m-m^\prime \right|},
\end{equation}

We now rewrite \eqref{eq:single_expectation_crosscorr} using $m^* = \min (m,m^\prime)$ and $\Delta m = | m - m^\prime|$
\begin{equation}
\E \left[ \VD_{p,\ell}^{n,m} \VD_{p,\ell}^{n,m^\prime} \right] = \left( 1 - \pncf^{m^*} \right) \left( 1 - \pncf^{M-m^*-\Delta m} \right) \pncf^{\Delta m}.
\end{equation}
If $\Delta m=0$ then $m=m^\prime$. For $\Delta m > 0$, there are two pairs of values $m, m^\prime$ for a given $m^*$, depending on whether $m$ or $m^\prime$ is the smaller of the two. The product $\VD_{p,\ell}^{n,m} \VD_{p,\ell}^{n,m^\prime}$ is either 0 or 1 so we work out the summation in \eqref{eq:cov_crosscorrelation} as

\begin{align}
\label{eq:double_sum_reordered}
\qquad\sum_{m=0}^{M{-}1} \! \sum_{m^\prime=0}^{M{-}1} \!\! \E \! \left[ \! \VD_{p,\ell}^{n,m} \VD_{p,\ell}^{n,m^\prime} \! \right] \! &= \!\!\!\underbrace{\sum_{m^*=1}^{M{-}1} \!\! \left(\! 1 {-} \pncf^{m^*} \!\right)\!\! \left(\! 1 {-} \pncf^{M{-}m^*} \!\right)}_{\Delta m = 0}\\
&\qquad+ \, 2 \!\!\! \sum_{\Delta m=1}^{M{-}2} \! \sum_{m^*=1}^{M{-}\Delta m{-}1} \!\! \left(\! 1 {-} \pncf^{m^*} \!\right) \!\! \left(\! 1 {-} \pncf^{M{-}m^*{-}\Delta m} \!\right) \! \pncf^{\Delta m}.
\end{align}

We simplify \eqref{eq:double_sum_reordered} by working out the sum over $m^*$
\begin{align}
\sum_{m^*=1}^{M-\Delta m-1}\!\!\!\! &\left( 1 - \pncf^{m^*} \right) \left( 1 - \pncf^{M-m^*-\Delta m} \right)\\
&= \sum_{m^*=1}^{M-\Delta m-1}\!\!\!\! \left( 1 - \pncf^{m^*} - \pncf^{M-m^*-\Delta m} + \pncf^{M-\Delta m} \right)\\
&= \left( 1 + \pncf^{M-\Delta m} \right) \left( M - \Delta m - 1\right) - \!\!\!\!\!\!\sum_{m^*=1}^{M-\Delta m-1}\!\!\!\!\!\! \pncf^{m^*}  - \pncf^{M-\Delta m}\!\!\!\!\sum_{m^*=1}^{M-\Delta m-1}\!\!\!\!\!\! \pncf^{-m^*}\\
&= \left( 1 + \pncf^{M-\Delta m} \right) \left( M - \Delta m - 1\right) - 2\frac{\pncf-\pncf^{M-\Delta m}}{1-\pncf}\\
&= M - 2\frac{\pncf}{\pcf} - 1 - \Delta m + \left( M + 2\pcf^{-1} - 1 \right) \pncf^{M-\Delta m} - \Delta m \pncf^{M-\Delta m}.
\end{align}
The double summation \eqref{eq:double_sum_reordered} now becomes
\begin{align}
\label{eq:double_sum_worked_out}
&\qquad \sum_{m=0}^{M{-}1} \sum_{m^\prime=0}^{M{-}1} \E \left[ \VD_{p,\ell}^{n,m} \VD_{p,\ell}^{n,m^\prime} \right] = M - 2\frac{\pncf}{\pcf} - 1 + \left( M + 2\pcf^{-1} - 1 \right) \pncf^{M}\\
&+\! 2\!\!\!\sum_{\Delta m=1}^{M-2}\!\!\! \left(\!\! \left(\!\! M\! - \! 2\frac{\pncf}{\pcf}\! - \! 1 \!\!\right)\! \pncf^{\Delta m}\!\! - \! \Delta m \pncf^{\Delta m}\! + \! \left(\! M + 2\pcf^{-1} \!\! - \! 1\! \right)\! \pncf^{M}\!\! - \!\Delta m \pncf^{M}\!\! \right)\!\!,
\end{align}
where the summation over $\Delta m$ can be worked out term by term
\begin{align}
&\!\!\!\!\!\!\!\!\!\!\!\!\!\!\!\!\!\!\!\!\!\!\!\!\!\!\!\!\!\!\!\!\!\!\!\!\sum_{\Delta m=1}^{M-2}\!\! \left(\!\! \left(\!\! M\! - \! 2\frac{\pncf}{\pcf}\! - \! 1 \!\!\right)\! \pncf^{\Delta m}\! - \Delta m \pncf^{\Delta m}\! + \! \left(\! M + 2\pcf^{-1} \! - \! 1\! \right)\! \pncf^{M}\! - \Delta m \pncf^{M}\! \right)\!\!\!\!\!\!\!\!\!\!\!\!\!\!\!\!\!\!\!\!\!\!\!\!\!\!\!\!\!\!\!\!\!\!\!\!\\
&\!\!\!\!\!\!\!\!\!\!\!\!\!\!\!\!\!\!\!\!\!\!\!\!\!\!\!\!\!\!\!\!\!\!\!\!\qquad\quad=\! \left(\!\! M\! - \! 2\frac{\pncf}{\pcf}\! - \! 1 \!\!\right)\!\frac{\pncf\!-\!\pncf^{M{-}1}}{1-\pncf} - \frac{(M{-}2)\pncf^{M}\!-\!(M{-}1)\pncf^{M{-}1}\!+\pncf}{(\pncf-1)^2}\!\!\!\!\!\!\!\!\!\!\!\!\!\!\!\!\!\!\!\!\!\!\!\!\!\!\!\!\!\!\!\!\!\!\!\!\\
&\!\!\!\!\!\!\!\!\!\!\!\!\!\!\!\!\!\!\!\!\!\!\!\!\!\!\!\!\!\!\!\!\!\!\!\!\qquad\qquad + \left( M-2 \right)\left( M + 2\pcf^{-1} - 1 \right) \pncf^{M} - (M-2)\frac{M - 1}{2} \pncf^M\!\!\!\!\!\!\!\!\!\!\!\!\!\!\!\!\!\!\!\!\!\!\!\!\!\!\!\!\!\!\!\!\!\!\!\!\\
&\!\!\!\!\!\!\!\!\!\!\!\!\!\!\!\!\!\!\!\!\!\!\!\!\!\!\!\!\!\!\!\!\!\!\!\!\qquad\quad=\! \frac{3\pncf^M \! - \! (M{+}1)\pncf^2 \! + \! (M{-}2)\pncf}{\pcf^2} \! + \! \left( M{-}2 \right)\!\!\left( \frac{M}{2}\! + 2\pcf^{-1}\! - \frac{1}{2} \right)\! \pncf^{M}.\!\!\!\!\!\!\!\!\!\!\!\!\!\!\!\!\!\!\!\!\!\!\!\!\!\!\!\!\!\!\!\!\!\!\!\!
\label{eq:deltam_sum}
\end{align}
Plugging \eqref{eq:deltam_sum} into \eqref{eq:double_sum_worked_out} gives us
\begin{align}
\sum_{m=0}^{M{-}1} \sum_{m^\prime=0}^{M{-}1} \E \left[ \VD_{p,\ell}^{n,m} \VD_{p,\ell}^{n,m^\prime} \right] &= 2\frac{3\pncf^M - (M+1)\pncf^2 + (M-2)\pncf}{\pcf^2}\\
+ M - 2\frac{\pncf}{\pcf} - &1 + \left(2\frac{M-2}{\pcf} + (M-1)\left( M + 2\pcf^{-1} - 1 \right) \right) \pncf^{M}
\end{align}

In \eqref{eq:variance_decomposed} (and as a consequence in \eqref{eq:total_variance}), \eqref{eq:cov_diffusioncorrelation} and \eqref{eq:cov_crosscorrelation} are the only terms that change compared to the corresponding expression derived for the term-by-term correlation, both appearing with a negative sign. This means that the variance of the difference of positions, and, as a consequence, any QoI's depending thereupon, is minimized by maximizing the sum of \eqref{eq:cov_diffusioncorrelation} and \eqref{eq:cov_crosscorrelation}. We rewrite this sum as
\begin{equation}
\label{eq:sum_of_roots}
C_1 \sqrt{\theta_\ell} + C_2 \sqrt{1-\theta_\ell}, \quad \text{with} \quad
C_1 = 2 \Delta t_{\ell-1}\sqrt{D_{\Delta t_\ell}D_{\Delta t_{\ell-1}}}
\end{equation}
\begin{equation}
\text{and} \quad C_2 = \Delta t_\ell \CV_{\Delta t_\ell} \sqrt{\frac{2 \Delta t_{\ell-1} D_{\Delta t_{\ell-1}}}{\V \left[ \sum_{m=0}^{M{-}1} \VD_{p,\ell}^{n,m^\prime} \right]}} \sum_{m=0}^{M{-}1} \sum_{m^\prime=0}^{M{-}1}   \E \left[ \VD_{p,\ell}^{n,m} \VD_{p,\ell}^{n,m^\prime} \right]
\end{equation}
positive and independent of $\theta_\ell$. To maximize \eqref{eq:sum_of_roots} we set the derivative in $\theta_\ell$ to zero
\begin{equation}
\label{eq:optimal_theta}
\frac{C_1}{2\sqrt{\theta_\ell}} - \frac{C_2}{2\sqrt{1-\theta_\ell}} = 0 \Leftrightarrow \frac{\theta_\ell}{1-\theta_\ell} = \frac{C_1^2}{C_2^2} \Leftrightarrow \theta_\ell = \frac{C_1^2}{C_1^2 + C_2^2}.
\end{equation}
Note that setting $\theta_\ell=1$ produces the term-by-term correlation from~\cite{Loevbak2021} as a special case of our new combined correlation. As we have shown that $C_1,C_2 > 0$, it is clear that the optimal value for $0 < \theta_\ell < 1$, and that the combined approach will produce pairs of particle paths which are more correlated than in the term-by-term approach, thus reducing the overall computational cost of the multilevel Monte Carlo simulation.

\section{Experiments}
\label{sec:experiments}

We now compare our new scheme with the term-by-term scheme. We first use the level strategy from in Section~\ref{sec:strategy}, observing similar performance from both schemes. We then show that the new scheme allows for a strategy with fewer levels, resulting in a large speedup compared to the term-by-term approach, where all levels are necessary. This speedup increases as $\epsilon$ decreases. The code for generating the results in this section can be found at \url{github.com/ELoevbak/APMLMC}.

\begin{table}
	\caption{Computing $\hat{Y}(t^*)$ with $\mathcal{M}(v) = \frac{1}{2}\left( \delta_{v,-1} \! + \delta_{v,1} \right)$ and term-by-term correlation. \label{tab:twospeedindependent}}
	\centering 
	\begin{tabular}{c | c c c c c | c}
		$\ell$ & $\Delta t_\ell$ & $P_\ell$ & \multicolumn{1}{c}{$\mathbb{E}[ \hat{F}_\ell \! - \! \hat{F}_{\ell-1} ]$} & $\mathbb{V}[ \hat{F}_\ell \! - \! \hat{F}_{\ell-1} ]$ & $\mathbb{V}[\hat{Y}_\ell]$ & $P_\ell C_\ell$ \\
		\hline
		0 & $5.0 \tabletimes 10^{-1}$ & $6.3 \tabletimes 10^8$ & $\phantom{-}9.90 \tabletimes 10^{-1}$ & $2.0 \tabletimes 10^{0\phantom{-}}$ & $3.1 \tabletimes 10^{-9}$ & $1.3 \tabletimes 10^7$\\
		1 & $1.0 \tabletimes 10^{-2}$ & $7.5 \tabletimes 10^7$ & $-1.25 \tabletimes 10^{-1}$ & $1.4 \tabletimes 10^{0\phantom{-}}$ & $1.9 \tabletimes 10^{-8}$ & $7.6 \tabletimes 10^7$\\
		2 & $5.0 \tabletimes 10^{-3}$ & $2.4 \tabletimes 10^7$ & $\phantom{-}1.08 \tabletimes 10^{-2}$ & $4.4 \tabletimes 10^{-1}$ & $1.8 \tabletimes 10^{-8}$ & $7.2 \tabletimes 10^7$\\
		3 & $2.5 \tabletimes 10^{-3}$ & $1.6 \tabletimes 10^7$ & $\phantom{-}2.83 \tabletimes 10^{-2}$ & $4.0 \tabletimes 10^{-1}$ & $2.5 \tabletimes 10^{-8}$ & $9.8 \tabletimes 10^7$\\
		4 & $1.3 \tabletimes 10^{-3}$ & $1.0 \tabletimes 10^7$ & $\phantom{-}2.91 \tabletimes 10^{-2}$ & $3.0 \tabletimes 10^{-1}$ & $3.0 \tabletimes 10^{-8}$ & $1.2 \tabletimes 10^8$\\
		5 & $6.3 \tabletimes 10^{-4}$ & $5.7 \tabletimes 10^6$ & $\phantom{-}2.09 \tabletimes 10^{-2}$ & $2.0 \tabletimes 10^{-1}$ & $3.4 \tabletimes 10^{-8}$ & $1.4 \tabletimes 10^8$\\
		6 & $3.1 \tabletimes 10^{-4}$ & $3.0 \tabletimes 10^6$ & $\phantom{-}1.26 \tabletimes 10^{-2}$ & $1.1 \tabletimes 10^{-1}$ & $3.6 \tabletimes 10^{-8}$ & $1.5 \tabletimes 10^8$\\
		7 & $1.6 \tabletimes 10^{-4}$ & $1.6 \tabletimes 10^6$ & $\phantom{-}6.87 \tabletimes 10^{-3}$ & $6.0 \tabletimes 10^{-2}$ & $3.8 \tabletimes 10^{-8}$ & $1.5 \tabletimes 10^8$\\
		8 & $7.8 \tabletimes 10^{-5}$ & $8.0 \tabletimes 10^5$ & $\phantom{-}3.46 \tabletimes 10^{-3}$ & $3.0 \tabletimes 10^{-2}$ & $3.8 \tabletimes 10^{-8}$ & $1.5 \tabletimes 10^8$\\
		9 & $3.9 \tabletimes 10^{-5}$ & $4.2 \tabletimes 10^5$ & $\phantom{-}1.60 \tabletimes 10^{-3}$ & $1.7 \tabletimes 10^{-2}$ & $4.0 \tabletimes 10^{-8}$ & $1.6 \tabletimes 10^8$\\
		10 & $1.9 \tabletimes 10^{-5}$ & $2.0 \tabletimes 10^5$ & $\phantom{-}1.06 \tabletimes 10^{-3}$ & $7.7 \tabletimes 10^{-3}$ & $3.9 \tabletimes 10^{-8}$ & $1.5 \tabletimes 10^8$\\
		11 & $9.8 \tabletimes 10^{-6}$ & $1.0 \tabletimes 10^5$ & $\phantom{-}5.94 \tabletimes 10^{-4}$ & $3.9 \tabletimes 10^{-3}$ & $3.9 \tabletimes 10^{-8}$ & $1.6 \tabletimes 10^8$\\
		12 & $4.9 \tabletimes 10^{-6}$ & $4.9 \tabletimes 10^4$ & $\phantom{-}3.40 \tabletimes 10^{-4}$ & $2.0 \tabletimes 10^{-3}$ & $4.1 \tabletimes 10^{-8}$ & $1.5 \tabletimes 10^8$\\
		13 & $2.4 \tabletimes 10^{-6}$ & $2.1 \tabletimes 10^4$ & $\phantom{-}4.55 \tabletimes 10^{-4}$ & $6.5 \tabletimes 10^{-4}$ & $3.1 \tabletimes 10^{-8}$ & $1.3 \tabletimes 10^8$\\
		14 & $1.2 \tabletimes 10^{-6}$ & $2.4 \tabletimes 10^4$ & $\phantom{-}2.07 \tabletimes 10^{-5}$ & $9.4 \tabletimes 10^{-4}$ & $4.0 \tabletimes 10^{-8}$ & $2.9 \tabletimes 10^8$\\
		15 & $6.1 \tabletimes 10^{-7}$ & $1.0 \tabletimes 10^3$ & $\phantom{-}2.91 \tabletimes 10^{-4}$ & $1.9 \tabletimes 10^{-4}$ & $1.9 \tabletimes 10^{-7}$ & $2.5 \tabletimes 10^7$\\
		\hline
		$\sum$ & \multicolumn{1}{c}{} & & $\phantom{-}9.81 \tabletimes 10^{-1}$ & & $6.6 \tabletimes 10^{-7}$ & $2.0 \tabletimes 10^9$
	\end{tabular}
\end{table}

\begin{table}
	\caption{Computing $\hat{Y}(t^*)$ with $\mathcal{M}(v) = \mathcal{N}(v;0,1)$ and term-by-term correlation. \label{tab:gaussianindependent}}
	\centering 
	\begin{tabular}{c | c c c c c | c}
		$\ell$ & $\Delta t_\ell$ & $P_\ell$ & \multicolumn{1}{c}{$\mathbb{E}[ \hat{F}_\ell \! - \! \hat{F}_{\ell-1} ]$} & $\mathbb{V}[ \hat{F}_\ell \! - \! \hat{F}_{\ell-1} ]$ & $\mathbb{V}[\hat{Y}_\ell]$ & $P_\ell C_\ell$ \\
		\hline
		0 & $5.0 \tabletimes 10^{-1}$ & $5.1 \tabletimes 10^8$ & $\phantom{-}9.90 \tabletimes 10^{-1}$ & $2.0 \tabletimes 10^{0\phantom{-}}$ & $3.9 \tabletimes 10^{-9}$ & $1.0 \tabletimes 10^7$\\
		1 & $1.0 \tabletimes 10^{-2}$ & $6.1 \tabletimes 10^7$ & $-1.25 \tabletimes 10^{-1}$ & $1.5 \tabletimes 10^{0\phantom{-}}$ & $2.4 \tabletimes 10^{-8}$ & $6.2 \tabletimes 10^7$\\
		2 & $5.0 \tabletimes 10^{-3}$ & $2.0 \tabletimes 10^7$ & $\phantom{-}1.08 \tabletimes 10^{-2}$ & $4.7 \tabletimes 10^{-1}$ & $2.3 \tabletimes 10^{-8}$ & $6.1 \tabletimes 10^7$\\
		3 & $2.5 \tabletimes 10^{-3}$ & $1.4 \tabletimes 10^7$ & $\phantom{-}2.82 \tabletimes 10^{-2}$ & $4.4 \tabletimes 10^{-1}$ & $3.2 \tabletimes 10^{-8}$ & $8.3 \tabletimes 10^7$\\
		4 & $1.3 \tabletimes 10^{-3}$ & $8.5 \tabletimes 10^6$ & $\phantom{-}2.88 \tabletimes 10^{-2}$ & $3.3 \tabletimes 10^{-1}$ & $3.9 \tabletimes 10^{-8}$ & $1.0 \tabletimes 10^8$\\
		5 & $6.3 \tabletimes 10^{-4}$ & $4.8 \tabletimes 10^6$ & $\phantom{-}2.12 \tabletimes 10^{-2}$ & $2.1 \tabletimes 10^{-1}$ & $4.4 \tabletimes 10^{-8}$ & $1.2 \tabletimes 10^8$\\
		6 & $3.1 \tabletimes 10^{-4}$ & $2.6 \tabletimes 10^6$ & $\phantom{-}1.23 \tabletimes 10^{-2}$ & $1.2 \tabletimes 10^{-1}$ & $4.8 \tabletimes 10^{-8}$ & $1.2 \tabletimes 10^8$\\
		7 & $1.6 \tabletimes 10^{-4}$ & $1.3 \tabletimes 10^6$ & $\phantom{-}6.57 \tabletimes 10^{-3}$ & $6.5 \tabletimes 10^{-2}$ & $4.9 \tabletimes 10^{-8}$ & $1.3 \tabletimes 10^8$\\
		8 & $7.8 \tabletimes 10^{-5}$ & $6.8 \tabletimes 10^5$ & $\phantom{-}4.12 \tabletimes 10^{-3}$ & $3.4 \tabletimes 10^{-2}$ & $5.0 \tabletimes 10^{-8}$ & $1.3 \tabletimes 10^8$\\
		9 & $3.9 \tabletimes 10^{-5}$ & $3.4 \tabletimes 10^5$ & $\phantom{-}2.11 \tabletimes 10^{-3}$ & $1.7 \tabletimes 10^{-2}$ & $5.0 \tabletimes 10^{-8}$ & $1.3 \tabletimes 10^8$\\
		10 & $1.9 \tabletimes 10^{-5}$ & $1.6 \tabletimes 10^5$ & $\phantom{-}1.20 \tabletimes 10^{-3}$ & $7.3 \tabletimes 10^{-3}$ & $4.6 \tabletimes 10^{-8}$ & $1.2 \tabletimes 10^8$\\
		11 & $9.8 \tabletimes 10^{-6}$ & $1.5 \tabletimes 10^5$ & $\phantom{-}6.58 \tabletimes 10^{-4}$ & $4.5 \tabletimes 10^{-3}$ & $2.9 \tabletimes 10^{-8}$ & $2.4 \tabletimes 10^8$\\
		12 & $4.9 \tabletimes 10^{-6}$ & $2.4 \tabletimes 10^4$ & $\phantom{-}5.23 \tabletimes 10^{-4}$ & $1.7 \tabletimes 10^{-3}$ & $7.1 \tabletimes 10^{-8}$ & $7.4 \tabletimes 10^7$\\
		13 & $2.4 \tabletimes 10^{-6}$ & $1.0 \tabletimes 10^3$ & $\phantom{-}1.11 \tabletimes 10^{-4}$ & $6.5 \tabletimes 10^{-4}$ & $6.5 \tabletimes 10^{-7}$ & $6.1 \tabletimes 10^6$\\
		\hline
		$\sum$ & \multicolumn{1}{c}{} & & $\phantom{-}9.82 \tabletimes 10^{-1}$ & & $1.2 \tabletimes 10^{-6}$ & $1.4 \tabletimes 10^9$
	\end{tabular}
\end{table}

\begin{table}
	\caption{Computing $\hat{Y}(t^*)$ with $\mathcal{M}(v) = \frac{1}{2}\left( \delta_{v,-1} \! + \delta_{v,1} \right)$ and combined correlation. \label{tab:twospeedcombined}}
	\centering 
	\begin{tabular}{c | c c c c c | c}
		$\ell$ & $\Delta t_\ell$ & $P_\ell$ & \multicolumn{1}{c}{$\mathbb{E}[ \hat{F}_\ell \! - \! \hat{F}_{\ell-1} ]$} & $\mathbb{V}[ \hat{F}_\ell \! - \! \hat{F}_{\ell-1} ]$ & $\mathbb{V}[\hat{Y}_\ell]$ & $P_\ell C_\ell$ \\
		\hline
		0 & $5.0 \tabletimes 10^{-1}$ & $4.5 \tabletimes 10^8$ & $\phantom{-}9.90 \tabletimes 10^{-1}$ & $2.0 \tabletimes 10^{0\phantom{-}}$ & $4.3 \tabletimes 10^{-9}$ & $9.1\tabletimes 10^6$\\
		1 & $1.0 \tabletimes 10^{-2}$ & $1.8 \tabletimes 10^7$ & $-1.25 \tabletimes 10^{-1}$ & $1.6 \tabletimes 10^{-1}$ & $8.8 \tabletimes 10^{-9}$ & $1.9 \tabletimes 10^7$\\
		2 & $5.0 \tabletimes 10^{-3}$ & $1.8 \tabletimes 10^7$ & $\phantom{-}1.04 \tabletimes 10^{-2}$ & $4.4 \tabletimes 10^{-1}$ & $2.5 \tabletimes 10^{-8}$ & $5.3 \tabletimes 10^7$\\
		3 & $2.5 \tabletimes 10^{-3}$ & $1.2 \tabletimes 10^7$ & $\phantom{-}2.85 \tabletimes 10^{-2}$ & $4.0 \tabletimes 10^{-1}$ & $3.4 \tabletimes 10^{-8}$ & $7.1 \tabletimes 10^7$\\
		4 & $1.3 \tabletimes 10^{-3}$ & $7.3 \tabletimes 10^6$ & $\phantom{-}2.86 \tabletimes 10^{-2}$ & $3.0 \tabletimes 10^{-1}$ & $4.2 \tabletimes 10^{-8}$ & $8.8 \tabletimes 10^7$\\
		5 & $6.3 \tabletimes 10^{-4}$ & $4.1 \tabletimes 10^6$ & $\phantom{-}2.06 \tabletimes 10^{-2}$ & $1.9 \tabletimes 10^{-1}$ & $4.7 \tabletimes 10^{-8}$ & $9.9 \tabletimes 10^7$\\
		6 & $3.1 \tabletimes 10^{-4}$ & $2.2 \tabletimes 10^6$ & $\phantom{-}1.24 \tabletimes 10^{-2}$ & $1.1 \tabletimes 10^{-1}$ & $5.0 \tabletimes 10^{-8}$ & $1.1 \tabletimes 10^8$\\
		7 & $1.6 \tabletimes 10^{-4}$ & $1.2 \tabletimes 10^6$ & $\phantom{-}6.84 \tabletimes 10^{-3}$ & $6.0 \tabletimes 10^{-2}$ & $5.2 \tabletimes 10^{-8}$ & $1.1 \tabletimes 10^8$\\
		8 & $7.8 \tabletimes 10^{-5}$ & $5.8 \tabletimes 10^5$ & $\phantom{-}3.55 \tabletimes 10^{-3}$ & $3.1 \tabletimes 10^{-2}$ & $5.3 \tabletimes 10^{-8}$ & $1.1 \tabletimes 10^8$\\
		9 & $3.9 \tabletimes 10^{-5}$ & $2.9 \tabletimes 10^5$ & $\phantom{-}1.38 \tabletimes 10^{-3}$ & $1.6 \tabletimes 10^{-2}$ & $5.4 \tabletimes 10^{-8}$ & $1.1 \tabletimes 10^8$\\
		10 & $1.9 \tabletimes 10^{-5}$ & $1.6 \tabletimes 10^5$ & $\phantom{-}9.31 \tabletimes 10^{-4}$ & $8.8 \tabletimes 10^{-3}$ & $5.6 \tabletimes 10^{-8}$ & $1.2 \tabletimes 10^8$\\
		11 & $9.8 \tabletimes 10^{-6}$ & $7.0 \tabletimes 10^4$ & $\phantom{-}6.67 \tabletimes 10^{-4}$ & $3.6 \tabletimes 10^{-3}$ & $5.2 \tabletimes 10^{-8}$ & $1.1 \tabletimes 10^8$\\
		12 & $4.9 \tabletimes 10^{-6}$ & $1.5 \tabletimes 10^4$ & $-1.13 \tabletimes 10^{-4}$ & $2.5 \tabletimes 10^{-3}$ & $1.6 \tabletimes 10^{-7}$ & $4.7 \tabletimes 10^7$\\
		13 & $2.4 \tabletimes 10^{-6}$ & $1.0 \tabletimes 10^3$ & $\phantom{-}2.24 \tabletimes 10^{-4}$ & $3.5 \tabletimes 10^{-4}$ & $3.5 \tabletimes 10^{-7}$ & $6.1 \tabletimes 10^6$\\
		\hline
		$\sum$ & \multicolumn{1}{c}{} & & $\phantom{-}9.79 \tabletimes 10^{-1}$ & & $1.0 \tabletimes 10^{-6}$ & $1.1 \tabletimes 10^9$
	\end{tabular}
\end{table}

\begin{table}
	\caption{Computing $\hat{Y}(t^*)$ with $\mathcal{M}(v) = \mathcal{N}(v;0,1)$ and combined correlation. \label{tab:gaussiancombined}}
	\centering 
	\begin{tabular}{c | c c c c c | c}
		$\ell$ & $\Delta t_\ell$ & $P_\ell$ & \multicolumn{1}{c}{$\mathbb{E}[ \hat{F}_\ell \! - \! \hat{F}_{\ell-1} ]$} & $\mathbb{V}[ \hat{F}_\ell \! - \! \hat{F}_{\ell-1} ]$ & $\mathbb{V}[\hat{Y}_\ell]$ & $P_\ell C_\ell$ \\
		\hline
		0 & $5.0 \tabletimes 10^{-1}$ & $5.1 \tabletimes 10^8$ & $\phantom{-}9.90 \tabletimes 10^{-1}$ & $2.0 \tabletimes 10^{0\phantom{-}}$ & $3.9 \tabletimes 10^{-9}$ & $1.0\tabletimes 10^7$\\
		1 & $1.0 \tabletimes 10^{-2}$ & $2.1 \tabletimes 10^7$ & $-1.25 \tabletimes 10^{-1}$ & $1.7 \tabletimes 10^{-1}$ & $8.0 \tabletimes 10^{-9}$ & $2.1 \tabletimes 10^7$\\
		2 & $5.0 \tabletimes 10^{-3}$ & $2.0 \tabletimes 10^7$ & $\phantom{-}1.04 \tabletimes 10^{-2}$ & $4.7 \tabletimes 10^{-1}$ & $2.3 \tabletimes 10^{-8}$ & $6.1 \tabletimes 10^7$\\
		3 & $2.5 \tabletimes 10^{-3}$ & $1.4 \tabletimes 10^7$ & $\phantom{-}2.84 \tabletimes 10^{-2}$ & $4.4 \tabletimes 10^{-1}$ & $3.2 \tabletimes 10^{-8}$ & $8.3 \tabletimes 10^7$\\
		4 & $1.3 \tabletimes 10^{-3}$ & $8.5 \tabletimes 10^6$ & $\phantom{-}2.88 \tabletimes 10^{-2}$ & $3.3 \tabletimes 10^{-1}$ & $3.9 \tabletimes 10^{-8}$ & $1.0 \tabletimes 10^8$\\
		5 & $6.3 \tabletimes 10^{-4}$ & $4.9 \tabletimes 10^6$ & $\phantom{-}2.08 \tabletimes 10^{-2}$ & $2.1 \tabletimes 10^{-1}$ & $4.4 \tabletimes 10^{-8}$ & $1.2 \tabletimes 10^8$\\
		6 & $3.1 \tabletimes 10^{-4}$ & $2.6 \tabletimes 10^6$ & $\phantom{-}1.22 \tabletimes 10^{-2}$ & $1.2 \tabletimes 10^{-1}$ & $4.7 \tabletimes 10^{-8}$ & $1.2 \tabletimes 10^8$\\
		7 & $1.6 \tabletimes 10^{-4}$ & $1.3 \tabletimes 10^6$ & $\phantom{-}6.87 \tabletimes 10^{-3}$ & $6.5 \tabletimes 10^{-2}$ & $4.9 \tabletimes 10^{-8}$ & $1.3 \tabletimes 10^8$\\
		8 & $7.8 \tabletimes 10^{-5}$ & $6.8 \tabletimes 10^5$ & $\phantom{-}3.65 \tabletimes 10^{-3}$ & $3.4 \tabletimes 10^{-2}$ & $5.0 \tabletimes 10^{-8}$ & $1.3 \tabletimes 10^8$\\
		9 & $3.9 \tabletimes 10^{-5}$ & $3.5 \tabletimes 10^5$ & $\phantom{-}2.22 \tabletimes 10^{-3}$ & $1.8 \tabletimes 10^{-2}$ & $5.2 \tabletimes 10^{-8}$ & $1.4 \tabletimes 10^8$\\
		10 & $1.9 \tabletimes 10^{-5}$ & $3.8 \tabletimes 10^5$ & $\phantom{-}1.03 \tabletimes 10^{-3}$ & $8.6 \tabletimes 10^{-3}$ & $2.3 \tabletimes 10^{-8}$ & $2.9 \tabletimes 10^8$\\
		11 & $9.8 \tabletimes 10^{-6}$ & $8.2 \tabletimes 10^4$ & $\phantom{-}4.10 \tabletimes 10^{-4}$ & $4.7 \tabletimes 10^{-3}$ & $5.8 \tabletimes 10^{-8}$ & $1.3 \tabletimes 10^8$\\
		12 & $4.9 \tabletimes 10^{-6}$ & $4.8 \tabletimes 10^4$ & $\phantom{-}3.27 \tabletimes 10^{-4}$ & $2.2 \tabletimes 10^{-3}$ & $4.8 \tabletimes 10^{-8}$ & $1.5 \tabletimes 10^8$\\
		13 & $2.4 \tabletimes 10^{-6}$ & $1.0 \tabletimes 10^3$ & $\phantom{-}3.91 \tabletimes 10^{-4}$ & $1.8 \tabletimes 10^{-4}$ & $1.8 \tabletimes 10^{-7}$ & $6.1 \tabletimes 10^6$\\
		\hline
		$\sum$ & \multicolumn{1}{c}{} & & $\phantom{-}9.81 \tabletimes 10^{-1}$ & & $6.6 \tabletimes 10^{-7}$ & $1.5 \tabletimes 10^9$
	\end{tabular}
\end{table}

As in~\cite{Loevbak2021}, we choose to estimate the squared particle displacement $F(x,v) = x^2$ at time $t^*=0.5$ using the Goldstein-Taylor model and Gaussian velocity distribution with $\CV=1$. In both cases we choose an initial condition $f(x,v,0) = \delta_{x,0} \mathcal{M}(v)$ and choose a desired root mean square error $E=0.001$. For each level $\ell$ we list the simulation time step $\Delta t_\ell$, number of samples $P_\ell$, estimated mean $\mathbb{E}[ \hat{F}_\ell \! - \! \hat{F}_{\ell-1} ]$ and variance $\mathbb{V}[ \hat{F}_\ell \! - \! \hat{F}_{\ell-1} ]$ of the difference samples, level variance $\mathbb{V}[\hat{Y}_\ell]$ and level cost $P_\ell C_\ell$ of simulations with the term-by-term scheme in Tables~\ref{tab:twospeedindependent}--\ref{tab:gaussianindependent}. We define $\hat{F}_{-1}\equiv0$, compute $\mathbb{V}[\hat{Y}_\ell] = P_\ell^{-1}\mathbb{V}[ \hat{F}_\ell \! - \! \hat{F}_{\ell-1} ]$ and specify the cost $C_\ell$ in terms of one single simulation with $\Delta t = 0.01$, i.e., the number of time steps performed divided by 50. The criteria for adding levels, distributing samples over the levels and stopping are the default ones as can be found in~\cite{Giles2015}. In Tables~\ref{tab:twospeedcombined}--\ref{tab:gaussiancombined} we list the same quantities for simulations using the new combined scheme, using the optimal value for $\theta_1$ as given by \eqref{eq:optimal_theta} and $\theta_\ell = 1$ for $\ell=2,\dots,L$. For the two speed model, we use the approach described in Section~\ref{sec:level0} for the level 0 simulation with $\Lambda=20$.

Our first observation from Tables~\ref{tab:twospeedindependent}--\ref{tab:gaussiancombined} is that the variance $\mathbb{V}[ \hat{F}_1 \! - \! \hat{F}_0 ]$ is much smaller when using the combined correlation approach. This means that fewer samples are needed at level 1 in Tables~\ref{tab:twospeedcombined}--\ref{tab:gaussiancombined} compared to Tables~\ref{tab:twospeedindependent}--\ref{tab:gaussianindependent}. The computed quantities from the different simulations are consistent to the order of magnitude of the requested $E$.

When comparing the total simulation cost for different correlation approaches, we do not observe a large difference between the tables, once the stochastic nature of the estimates in the stop-criterion is taken into account. The main reason for this is that the cost savings at level 1 are negligible as part of the total simulation cost. At this point it should be noted that the variance of the samples at level 1 in Tables~\ref{tab:twospeedcombined}--\ref{tab:gaussiancombined} is noticeably smaller than that at the subsequent levels. Given the base assumption for using MLMC that the variance should decrease with increasing $\ell$~\cite{Giles2008}, this indicates that the given sequence of levels is not optimal for the combined correlation approach.

\setlength{\figurewidth}{0.9\textwidth}
\setlength{\figureheight}{0.35\figurewidth}

In~\cite{Loevbak2020}, preliminary results indicated that leaving out a number of intermediate levels, when using a variant of the combined correlation, results in a significant speedup of the overall computational cost of the multilevel simulation. The same speedup was not observed for the term-by-term approach. We perform a more rigorous comparison of the effect of different level selection strategies for the new combined correlation and the term-by-term correlation from~\cite{Loevbak2021}. To this end we perform the simulations in Tables \ref{tab:twospeedindependent}--\ref{tab:gaussiancombined} for different values of $\epsilon$. For each $\epsilon$ we vary $\Delta t_1$, while keeping $\Delta t_0=0.5$ and setting $\Delta t_{\ell-1} = 2 \Delta t_\ell$ for $\ell > 1$.

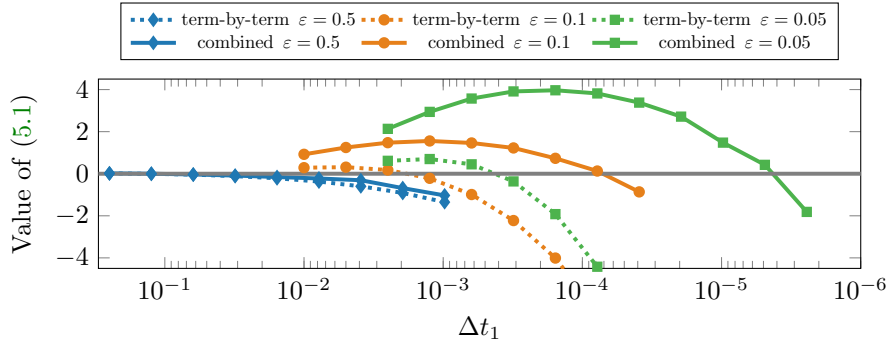
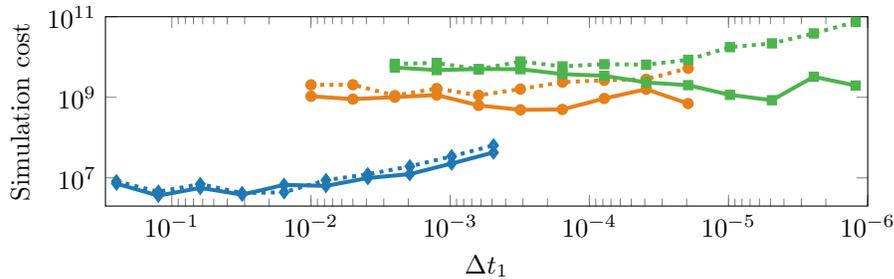
\begin{figure}
	\begin{subfigure}{\linewidth}
		\centering
		\begin{minipage}{0.95\linewidth}
			\flushright
			\begin{tikzpicture}

\begin{axis}[
xmin=1e-06, xmax=0.3,
ymin=-4.5, ymax=4.5,
xmode=log,
x dir=reverse,
xlabel=$\Delta t_1$,
ylabel={Value of \eqref{eq:leaveoutlevel}},
width=\figurewidth,
height=\figureheight,
legend style={nodes={scale=0.7, transform shape}, at={(0.5,1.1)},anchor=south},
legend columns=3,
x grid style={white!69.01960784313725!white},
y grid style={white!69.01960784313725!white},
every axis plot/.append style={line width=1.5pt}
]

\addplot[color0, dotted, mark=diamond*, mark size=1.6, mark options={solid}]
table {%
0.25	         0.015698818
0.125    	-0.003449317
0.0625    	-0.042845124
0.03125	    -0.114926393
0.015625  	-0.216658521
0.0078125	-0.371762364
0.00390625	-0.59856141
0.001953125	-0.912385194
0.000976563	-1.349462318
};

\addplot[color1, dotted, mark=*, mark size=1.4, mark options={solid}]
table {%
0.01	    0.287488861
0.005	0.31198284
0.0025	0.176633777
0.00125	    -0.217215865
0.000625  	-0.988220711
0.0003125	-2.225467365
0.00015625	-4.008553789
0.000078125	-6.637082196
3.90625E-05	-10.33393429
};

\addplot[color2, dotted, mark=square*, mark size=1.2, mark options={solid}]
table {%
0.0025	     0.614460079
0.00125	     0.698502363
0.000625	     0.447863763
0.0003125	-0.36490724
0.00015625	-1.919859814
0.000078125	-4.420668355
3.90625E-05	-8.106205212
1.95313E-05	-13.58799258
9.76563E-06	-21.03107315
4.88281E-06	-32.18112878
2.44141E-06	-47.45308544
};

\addplot[color0, mark=diamond*, mark size=1.6, mark options={solid}]
table {%
0.25 	0.018139723
0.125	0.000801407
0.0625	 -0.035524572
0.03125	 -0.087909742
0.015625	 -0.156512978
0.0078125	-0.225495327
0.00390625	-0.308929681
0.001953125	-0.683296389
0.000976563	-1.031122076
};

\addplot[color1, mark=*, mark size=1.4, mark options={solid}]
table {%
0.01					0.923109900580663
0.005					1.24736545090525
0.0025					1.47623829461663
0.00125					1.55786486041444
0.000625					1.462516258232
0.0003125					1.22584660697881
0.00015625					0.732873073334852
7.8125E-05					0.128729512906755
3.90625E-05					-0.863030531425052
};

\addplot[color2, mark=square*, mark size=1.2, mark options={solid}]
table {%
0.0025	2.133927644
0.00125	2.940418923
0.000625	3.574155804
0.0003125	3.912910315
0.00015625	3.967668622
0.000078125	3.81743605
3.90625E-05	3.380658288
1.95313E-05	2.714490813
9.76563E-06	1.478714972
4.88281E-06	0.422809319
2.44141E-06	-1.814253567
};

\addplot[color=gray] coordinates {(1e-06,0) (0.3,0)};

\legend{term-by-term\, $\epsilon = 0.5$, term-by-term\, $\epsilon = 0.1$, term-by-term\, $\epsilon = 0.05$, combined\, $\epsilon = 0.5$, combined\, $\epsilon = 0.1$, combined\, $\epsilon = 0.05$}

\end{axis}

\end{tikzpicture}
		\end{minipage}
		\caption{Computing the value of $\Delta t_1$ which sets \eqref{eq:leaveoutlevel} to zero using simulation estimated variances.\label{fig:GTcosttheoretical}}
	\end{subfigure}\newline
	\begin{subfigure}{\linewidth}
		\centering
		\begin{minipage}{0.95\linewidth}
			\flushright
			\begin{tikzpicture}

\begin{axis}[
xmin=1e-06, xmax=0.3,
ymin=2e6, ymax=1e11,
xmode=log,
ymode=log,
x dir=reverse,
xlabel=$\Delta t_1$,
ylabel={Simulation cost},
width=\figurewidth,
height=\figureheight,
legend style={nodes={scale=0.65, transform shape}},
legend pos=north west,
x grid style={white!69.01960784313725!white},
y grid style={white!69.01960784313725!white},
every axis plot/.append style={line width=1.5pt}
]

\addplot[color0, dotted, mark=diamond, mark size=1.6, mark options={solid}]
table {%
2.50E-01	 8.22E+06
1.25E-01	 4.53E+06
6.25E-02	 7.00E+06
3.13E-02	 4.19E+06
1.56E-02	 4.40E+06
7.81E-03	 8.78E+06
3.91E-03	 1.23E+07
1.95E-03	 1.93E+07
9.77E-04	 3.41E+07
4.88E-04	 6.42E+07
};

\addplot[color1, dotted, mark=*, mark size=1.4, mark options={solid}]
table {%
1.00E-02	 2.04E+09
5.00E-03	 2.05E+09
2.50E-03	 1.10E+09
1.25E-03	 1.65E+09
6.25E-04	 1.12E+09
3.13E-04	 1.60E+09
1.56E-04	 2.38E+09
7.81E-05	 2.65E+09
3.91E-05	 2.80E+09
1.95E-05	5.26E+09
};

\addplot[color2, dotted, mark=square*, mark size=1.2, mark options={solid}]
table {%
2.50E-03	 6.77E+09
1.25E-03	 7.13E+09
6.25E-04	 4.93E+09
3.13E-04	 7.74E+09
1.56E-04	 5.86E+09
7.81E-05	 6.63E+09
3.91E-05	 6.45E+09
1.95E-05	8.52E+09
9.77E-06	1.76E+10
4.88E-06	2.18E+10
2.44E-06	3.85E+10
1.22E-06	7.47E+10
};

\addplot[color0, mark=diamond, mark size=1.6, mark options={solid}]
table {%
0.25	     7.21E+06
0.125	 3.63E+06
0.0625	 5.65E+06
0.03125	 3.84E+06
0.015625	 6.71E+06
0.0078125	6.33E+06
0.00390625	9.90E+06
0.001953125	1.24E+07
0.000976563	2.27E+07
0.000488281	4.27E+07
};

\addplot[color1, mark=*, mark size=1.4, mark options={solid}]
table {%
0.01	     1.06E+09
0.005	 9.03E+08
0.0025	 1.01E+09
0.00125	 1.14E+09
0.000625 	6.29E+08
0.0003125	4.89E+08
0.00015625	4.99E+08
0.000078125	9.35E+08
3.90625E-05	1.59E+09
1.95313E-05	6.99E+08
};

\addplot[color2, mark=square*, mark size=1.2, mark options={solid}]
table {%
0.0025	5.49E+09
0.00125	4.74E+09
0.000625	5.02E+09
0.0003125	4.98E+09
0.00015625	3.77E+09
0.000078125	3.42E+09
3.90625E-05	2.34E+09
1.95313E-05	2.00E+09
9.76563E-06	1.15E+09
4.88281E-06	8.47E+08
2.44141E-06	3.24E+09
1.2207E-06	1.96E+09
};


\end{axis}

\end{tikzpicture}
		\end{minipage}
		\caption{Total simulation cost in terms of samples with $\Delta t = 0.01$.\label{fig:GTcostactual}}
	\end{subfigure}
	\caption{Comparing different level strategies for the two-speed model by keeping $\Delta t_0$ fixed and varying $\Delta t_1$. Figure~\ref{fig:GTcosttheoretical} consistently shows an optimal $\Delta t_1$ some orders of magnitude smaller than $\epsilon^2$. Figure~\ref{fig:GTcostactual} also shows that the combined correlation consistently outperforms the term-by-term correlation. \label{fig:GTcost}}
\end{figure}

\begin{figure}
	\begin{subfigure}{\linewidth}
		\centering
		\begin{minipage}{0.95\linewidth}
			\flushright
			\begin{tikzpicture}

\begin{axis}[
xmin=1e-06, xmax=0.3,
ymin=-4.5, ymax=4.5,
xmode=log,
x dir=reverse,
xlabel=$\Delta t_1$,
ylabel={Value of \eqref{eq:leaveoutlevel}},
width=\figurewidth,
height=\figureheight,
legend style={nodes={scale=0.7, transform shape}, at={(0.5,1.1)},anchor=south},
legend columns=3,
x grid style={white!69.01960784313725!white},
y grid style={white!69.01960784313725!white},
every axis plot/.append style={line width=1.5pt}
]

\addplot[color0, dotted, mark=diamond*, mark size=1.6, mark options={solid}]
table {%
0.25	     0.018671938
0.125	 0.001515119
0.0625	 -0.040929363
0.03125	 -0.116045254
0.015625	 -0.23991171
0.0078125	-0.423146981
0.00390625	-0.688634318
0.001953125	-1.059533475
0.000976563	-1.617532137
};

\addplot[color1, dotted, mark=*, mark size=1.4, mark options={solid}]
table {%
0.01 	0.306741904
0.005	 0.340877847
0.0025	 0.210991317
0.00125	 -0.189546805
0.000625 	-0.989699569
0.0003125	-2.218810783
0.00015625	-4.095579791
0.000078125	-6.796308744
3.90625E-05	-10.69965044
};

\addplot[color2, dotted, mark=square*, mark size=1.2, mark options={solid}]
table {%
0.0025	    0.62489029
0.00125	    0.713703439
0.000625 	0.462527024
0.0003125	-0.351606698
0.00015625	-1.920151466
0.000078125	-4.444395221
3.90625E-05	-8.18262805
1.95313E-05	-13.62957984
9.76563E-06	-21.32579127
4.88281E-06	-32.1435614
2.44141E-06	-47.62691558
};

\addplot[color0, mark=diamond*, mark size=1.6, mark options={solid}]
table {%
0.25   	0.020780055
0.125	0.004265834
0.0625	 -0.033441675
0.03125	 -0.104059967
0.015625 	-0.221333963
0.0078125	-0.394932885
0.00390625	-0.629220246
0.001953125	-1.034544459
0.000976563	-1.531189971
};

\addplot[color1, mark=*, mark size=1.4, mark options={solid}]
table {%
0.01					0.957425325702313
0.005					1.29926604787839
0.0025					1.54988350626228
0.00125					1.64109400849911
0.000625					1.55333932312963
0.0003125					1.30410998997723
0.00015625					0.790195341340787
7.8125E-05					0.038175609802238
3.90625E-05					-0.913413255872072
};

\addplot[color2, mark=square*, mark size=1.2, mark options={solid}]
table {%
0.0025	 2.152538096
0.00125 	2.979816755
0.000625 	3.618124082
0.0003125	3.966296239
0.00015625	 4.018916885
0.000078125 	3.821539783
3.90625E-05	 3.465696017
1.95313E-05	    2.698007635
9.76563E-06 	1.883720358
4.88281E-06 	0.123845771
2.44141E-06 	-1.930519455
};

\addplot[color=gray] coordinates {(1e-06,0) (0.3,0)};

\legend{term-by-term\, $\epsilon = 0.5$, term-by-term\, $\epsilon = 0.1$, term-by-term\, $\epsilon = 0.05$, combined\, $\epsilon = 0.5$, combined\, $\epsilon = 0.1$, combined\, $\epsilon = 0.05$}

\end{axis}

\end{tikzpicture}
		\end{minipage}
		\caption{Computing the value of $\Delta t_1$ which sets \eqref{eq:leaveoutlevel} to zero using simulation estimated variances.\label{fig:gausscosttheoretical}}
	\end{subfigure}\newline
	\begin{subfigure}{\linewidth}
		\centering
		\begin{minipage}{0.95\linewidth}
			\flushright
			\begin{tikzpicture}

\begin{axis}[
xmin=1e-06, xmax=0.3,
ymin=2e6, ymax=1e11,
xmode=log,
ymode=log,
x dir=reverse,
xlabel=$\Delta t_1$,
ylabel={Simulation cost},
width=\figurewidth,
height=\figureheight,
legend style={nodes={scale=0.65, transform shape}},
legend pos=north west,
x grid style={white!69.01960784313725!white},
y grid style={white!69.01960784313725!white},
every axis plot/.append style={line width=1.5pt}
]

\addplot[color0, dotted, mark=diamond, mark size=1.6, mark options={solid}]
table {%
0.25   	1.23E+07
0.125	1.08E+07
0.0625	5.37E+06
0.03125	5.80E+06
0.015625  	1.00E+07
0.0078125	1.60E+07
0.00390625	1.63E+07
0.001953125	3.52E+07
0.000976563	4.72E+07
0.000488281	1.02E+08
};

\addplot[color1, dotted, mark=*, mark size=1.4, mark options={solid}]
table {%
0.01	        1.48E+09
0.005	    3.60E+09
0.0025   	1.39E+09
0.00125	    1.81E+09
0.000625	    1.23E+09
0.0003125	1.37E+09
0.00015625	1.87E+09
0.000078125	2.99E+09
3.90625E-05	6.82E+09
1.95313E-05	5.80E+09
};

\addplot[color2, dotted, mark=square*, mark size=1.2, mark options={solid}]
table {%
0.0025    	7.65E+09
0.00125	    6.30E+09
0.000625	    5.02E+09
0.0003125	6.93E+09
0.00015625	9.11E+09
0.000078125	6.77E+09
3.90625E-05	7.29E+09
1.95313E-05	8.31E+09
9.76563E-06	1.31E+10
4.88281E-06	2.23E+10
2.44141E-06	3.84E+10
1.2207E-06	7.48E+10
};

\addplot[color0, mark=diamond, mark size=1.6, mark options={solid}]
table {%
0.25	    2.73E+07
0.125	8.91E+06
0.0625	7.42E+06
0.03125	9.05E+06
0.015625	    8.94E+06
0.0078125	1.95E+07
0.00390625	1.35E+07
0.001953125	2.49E+07
0.000976563	4.45E+07
0.000488281	8.65E+07
};

\addplot[color1, mark=*, mark size=1.4, mark options={solid}]
table {%
0.01								1.48E+09
0.005								1.07E+09
0.0025								1.10E+09
0.00125								1.25E+09
0.000625								9.48E+08
0.0003125								8.53E+08
0.00015625								7.14E+08
7.8125E-05								6.55E+08
3.90625E-05								6.51E+08
1.953125E-05								6.62E+08
};

\addplot[color2, mark=square*, mark size=1.2, mark options={solid}]
table {%
0.0025	6.05E+09
0.00125	6.11E+09
0.000625	4.54E+09
0.0003125	4.44E+09
0.00015625	5.13E+09
0.000078125	3.17E+09
3.90625E-05	3.04E+09
1.95313E-05	4.70E+09
9.76563E-06	2.50E+09
4.88281E-06	1.37E+09
2.44141E-06	2.36E+09
1.2207E-06	2.73E+09
};

\end{axis}

\end{tikzpicture}
		\end{minipage}
		\caption{Total simulation cost in terms of samples with $\Delta t = 0.01$.\label{fig:gausscostactual}}
	\end{subfigure}
	\caption{Comparing different level strategies with Gaussian velocities by keeping $\Delta t_0$ fixed and varying $\Delta t_1$. Figure~\ref{fig:gausscosttheoretical} consistently shows an optimal $\Delta t_1$ some orders of magnitude smaller than $\epsilon^2$. Figure~\ref{fig:gausscostactual} shows that the combined correlation consistently outperforms the term-by-term correlation. \label{fig:gausscost}}
\end{figure}
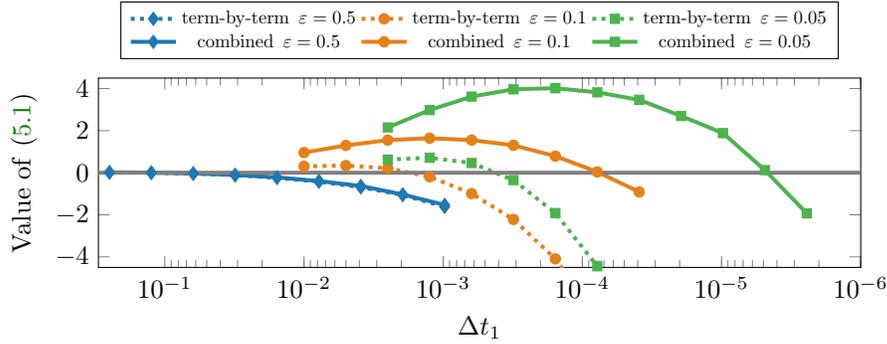
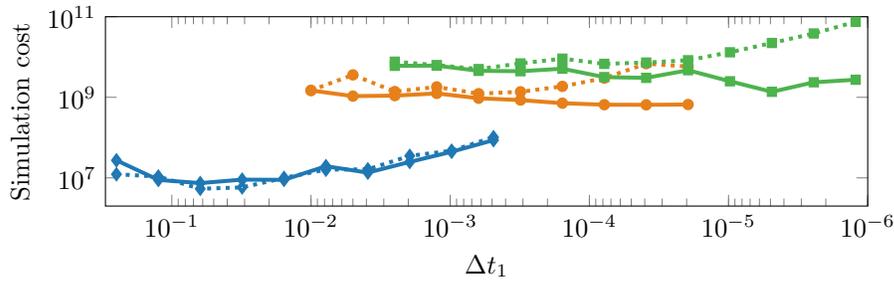

We consider two metrics to determine how good a given sequence of levels is. First off, we consider the total simulation cost of each sequence of levels, i.e., the sum $\sum_{\ell=0}^L P_\ell C_\ell$. This cost is computed using the actual number of samples $P_\ell$ and the analytically computed cost per sample $C_\ell$ at each level. Secondly, we make use of a theoretical result stating that it is beneficial to leave out a level $\ell$ if~\cite{Giles2015}
\begin{equation}
\label{eq:leaveoutlevel}
\sqrt{\mathbb{V}[ \hat{F}_\ell \! - \! \hat{F}_{\ell-1}]C_\ell} + \sqrt{\mathbb{V}[ \hat{F}_{\ell+1} \! - \! \hat{F}_{\ell}]C_{\ell+1}} - \sqrt{\mathbb{V}[ \hat{F}_{\ell+1} \! - \! \hat{F}_{\ell-1}]C^\prime_{\ell+1}}
\end{equation}
is positive, with $C^\prime_{\ell+1}$ the cost of performing a correlated simulation with time steps $\Delta t_{\ell-1}$ and $\Delta t_{\ell+1}$. To evaluate \eqref{eq:leaveoutlevel}, we use the variance estimates computed by the simulations. Both quantities are shown in Figures~\ref{fig:GTcost}--\ref{fig:gausscost} for a range of values $\Delta t_1$.

As Figures~\ref{fig:GTcost}--\ref{fig:gausscost} are similar, we draw the same two conclusions for both figures:
\begin{enumerate}
	\item We see in sub-figure a that the combined correlation works optimally with a much smaller level 1 time step size than the term-by-term correlation. We thus require fewer levels to bridge the gap between diffusive simulations with large time steps and kinetic simulations with time steps that are fine enough to meet the given error tolerance. We also see that the gap between the optimal level strategies for both correlations increases as $\epsilon$ decreases.
	\item We see in sub-figure b that the combined correlation consistently produces results at lower total simulation cost than the term-by-term correlation. We also notice that this cost is highly variable. This variation is due to the finer levels, where fewer samples are used. This reduced number of samples results in variation in the variance estimates which in turn leads to variation in the estimates for the number of samples. While this cost will likely be minimal for the optimal level strategy it does not appear to be highly sensitive to the selected value for $\Delta t_1$.  Simply setting $\Delta t_1$ to be some orders of magnitude smaller than $\epsilon^2$ suffices when selecting a level hierarchy.
\end{enumerate}
  
\begin{table}
\caption{Mean computed QoI per $\epsilon$ (with standard deviation between brackets) and speedup of best result compared with best term-by-term result.\label{tab:computedresults}}
\centering
\begin{tabular}{c | c | c | c | c|}
	&  \multicolumn{2}{c|}{$\mathcal{M}(v) = \frac{1}{2}\left( \delta_{v,-1} \! + \delta_{v,1} \right)$}  &  \multicolumn{2}{c|}{$\mathcal{M}(v) = \mathcal{N}\left( v ; 0,1 \right)$} \\
$\epsilon$ & Result & Speedup & Result & Speedup \\
\hline
0.5 & $5.68 \tabletimes 10^{-1}$ ($6.4 \tabletimes 10^{-4}$) & 1.2 & $5.68 \tabletimes 10^{-1}$ ($5.9 \tabletimes 10^{-4}$) & 0.7\\
0.1 & $9.79 \tabletimes 10^{-1}$ ($7.8 \tabletimes 10^{-4}$) & 2.2 & $9.80 \tabletimes 10^{-1}$ ($7.0 \tabletimes 10^{-4}$) & 1.9\\
0.05 & $9.95 \tabletimes 10^{-1}$ ($6.4 \tabletimes 10^{-4}$) & 5.8 & $9.95 \tabletimes 10^{-1}$ ($8.3 \tabletimes 10^{-4}$) & 3.7\\
\end{tabular}
\end{table}

In Table \ref{tab:computedresults} we present the observed speedup when comparing the most efficient computations achieved by both schemes. We observe that the speedup increases with decreasing epsilon as one would expect when approaching the diffusive limit. In the supplementary materials we include tables in the same form as Tables \ref{tab:twospeedindependent}--\ref{tab:gaussiancombined} for each simulation in Figures \ref{fig:GTcost}--\ref{fig:gausscost}. To verify the two speed simulation consistency, we also report the mean and standard deviation of the computed QoI averaged over each combined correlation curve in the figures. From this table, we conclude that the results computed with both distributions are the same, down to the specified RMSE $E$. This indicates that any bias caused by the run-length distribution sampling approach described in Section~\ref{sec:runlengthsampling} is not observable in the computed results. 

In conclusion, we propose to modify the level sequence proposed in~\cite{Loevbak2021} by reducing $\Delta t_1 \ll \epsilon^2$, so that a large number of unnecessary levels are skipped. A good rule of thumb seems to be one or two orders of magnitude smaller than $\epsilon^2$, but further work is needed in order to check whether this applies outside of the tested range for $\epsilon$. Apart from stating that $\Delta t_0 \gg \epsilon^2 \gg \Delta t_1$, we thus conclude that the primary factors in determining a sequence of levels for a given simulation will likely be the constraints given by the problem. On the one hand, non-homogeneous simulations may constrain $\Delta t_0$ in order to capture parameter dependence. On the other hand, practical applications often have a very high simulation cost and, as a consequence, a large tolerance $E$, which will result in a larger $\Delta t_1$.

\section{Conclusions}
\label{sec:conclusions}

We have demonstrated that existing APMLMC schemes can be improved by taking into account that coarse simulation diffusion substitutes fine simulation transport-collision dynamics. This requires the summation of fine simulation velocities to generate correlated coarse simulation diffusion. We have described in detail how to achieve this correlation in practice for arbitrary symmetric velocity distributions. When the velocity distribution is Gaussian, the proposed approach introduces no bias compared to simulating only at the finest level, while significantly accelerating computations. For other velocity distributions, a bias is present which can be significantly decreased by modifying the level 0 simulation. This modification comes at additional cost, however, which increases as the required bias decreases.

While the presented scheme is promising, further work is needed to adapt it to practical use cases. An extension to non-symmetric velocity distributions will be needed to cover cases where net drift effects are present in the model. The modified simulation at level 0, while performing less computation than simulating with a time step $\Delta t_1$, still performs a non-negligible amount of computation which is not taken into account in the considered cost-analysis. A large amount of memory is also required in its implementation. To analyze the significance these effects, we intend to perform a more efficient implementation in a compiled language in future work. Another potential area for future work, is to apply the approach in Section~\ref{sec:level0} to the scheme presented in~\cite{Mortier2020} and compare both schemes.

\appendix

\section{Sampling the run length statistics}
\label{app:run_length}

\renewcommand{\pncf}{{p_{nc}}}
\renewcommand{\pcf}{{p_c}}

In this appendix we briefly describe the approach presented in~\cite{Fu1994} for generating the probability distributions of $E_{M-1,\lambda^*{-}1}$ and $G_{M-1,\lambda^*{-}1}$, while omitting the derivation of the method. For further details we refer to the original publication. It is important to note that we have made significant modifications to the notation used in~\cite{Fu1994} to maintain the internal consistency of this paper. We have also opted to transpose all matrices and vectors so as to represent vectors as columns.

Given $M-1$ potential collisions with collision probability $\pcf$ and no-collision probability $\pncf$, $E_{M-1,\lambda^*{-}1}$ counts the number of times where precisely $\lambda^*{-}1$ such moments occur in sequence without a collision and $G_{M-1,\lambda^*{-}1}$ counts the number of times where at least $\lambda^*{-}1$ such moments occur in sequence without a collision. A similar approach can be used to calculate both statistics. First an $r$-dimensional state space $\Omega$ is defined, with $r$ a finite integer value depending on $M$, $\lambda^*$ and the statistic being generated. On this state space, a finite Markov chain is defined, with an $r \times r$ transition matrix $A$. Given a suitable initial state $\pi_0 \in \R^r$ and a vector $U \in \{0,1\}^r$, one can then compute the probability of the statistic $R_{M-1,\lambda^*{-}1} = E_{M-1,\lambda^*{-}1} \vee G_{M-1,\lambda^*{-}1}$ taking a given value $\omega = 0,\dots, l=\left \lfloor \frac{M}{\lambda^*} \right \rfloor$ as
\begin{equation}
\label{eq:run_probability}
P(R_{M-1,\lambda^*{-}1}=\omega) = U^\prime_\omega A^{M{-}1} \pi_0.
\end{equation}
For $E_{M-1,\lambda^*{-}1}$, set $r = (l+1)(\lambda^*+2)-1$ and the $r \times r$ transition matrix $A$ of the Markov chain takes the form

\newcommand{\mcdots}{\cdots}
\newcommand{\mddots}{\ddots}
\newcommand{\mvdots}{\vdots}
\begin{equation}
\def\arraystretch{0.8}
A = \left[\arraycolsep=1pt
    \begin{array}{c ;{2pt/2pt} c c c ;{2pt/2pt} c c ;{2pt/2pt} c c c ;{2pt/2pt} c c ;{2pt/2pt} c c c}
        \overmat{1}{\;\pncf\;} & \overmat{\lambda^*}{\hphantom{\pncf} & \;\hphantom{\ddots} & \hphantom{\pncf}} & \overmat{2}{\pncf & \hphantom{\pncf}} & \overmat{\lambda^*}{\hphantom{\pncf} & \;\hphantom{\ddots} & \hphantom{\pncf}} & \overmat{2}{\hphantom{\pncf} & \hphantom{\pncf}} & \overmat{\lambda^*}{\hphantom{\pncf} & \;\hphantom{\ddots} & \hphantom{\pncf}} \\
        \hdashline[2pt/2pt]
        \pcf & \pcf & \mcdots & \pcf & & & & & & & & \\
        & \pncf & & & & & & & & & & \\
        & & \!\mddots & & & & & & & & & \\
        \hdashline[2pt/2pt]
        & & & \pncf & 0 & & & & & & &\\
        & & & & & \pncf & & & & \pncf & & & \\
        \hdashline[2pt/2pt]
		& & & & \pcf & \pcf & \pcf & \mcdots & \pcf & & & & & \!\mddots \\
		& & & & & & \pncf & & & & & & \!\mddots \\
		& & & & & & & \!\mddots & & & & \!\mddots \\
		\hdashline[2pt/2pt]
		& & & & & & & & \pncf & & \!\mddots & & & \\
		& & & & & & & & & \!\mddots & & \\
		\hdashline[2pt/2pt]
		& & & & & & & & \!\mddots &\pcf & \pcf & \pcf & \mcdots & 0 \\		
		& & & & & & & \!\mddots & & & & \pncf & & \\
		& & & & & & \!\mddots & & & & & & \!\mddots & 1 \\		
    \end{array}
\right].
\end{equation}

The initial state $\pi_0$ and observation vector $u_\omega$ for a given $\omega$ are given by 
\begin{equation}
\pi_0 = \begin{bmatrix}
0 \\
1 \\
0 \\
\mvdots \\
0
\end{bmatrix}, \quad
U^\prime_0 = \begin{bmatrix}
\overmat{\lambda^*+1}{1 & \mcdots & 1} & \overmat{(\lambda^*+2)l}{0 & \mcdots & 0}
\end{bmatrix} \quad \text{and}
\end{equation}
\vspace{0pt}
\begin{equation}
U^\prime_\omega = \begin{bmatrix}
\overmat{(\lambda^*+2)\omega-1}{0 & \mcdots & 0} & \overmat{\lambda^*+2}{1 & \mcdots & 1} & \overmat{(\lambda^*+2)(l-\omega)}{0 & \mcdots & 0}
\end{bmatrix} , \quad \omega = 1,\dots l.
\end{equation}
For $G_{M-1,\lambda^*-1}$, set $r = (l+1)(\lambda^*+1)-1$ and the transition matrix $A$ takes the form

\begin{equation}
A = \left[\arraycolsep=2pt
    \begin{array}{c c c ;{2pt/2pt} c ;{2pt/2pt} c c c ;{2pt/2pt} c ;{2pt/2pt} c c c c}
        \overmat{\lambda^*}{\;\;\pcf & \;\mcdots & \pcf\;\;} & \overmat{1}{\hphantom{\;\mddots\;\;}} & \overmat{\lambda^*}{\hphantom{\pncf} & \;\hphantom{\mddots} & \hphantom{\pncf}} & \overmat{1}{\;\hphantom{\pncf}\;} & \overmat{\lambda^*}{\hphantom{\pncf} & \;\hphantom{\mddots} & \hphantom{\pncf}} \\
        \pncf & & & & & & & & & & \\
        & \!\mddots & & & & & & & & & \\
        \hdashline[2pt/2pt]
        & & \pncf & \pncf & & & & & & &\\
        \hdashline[2pt/2pt]
		& & & \pcf & \pcf & \mcdots & \pcf & & & & \!\mddots \\
		& & & & \pncf & & & & & \!\mddots \\
		& & & & & \!\mddots & & & \!\mddots \\
		\hdashline[2pt/2pt]
		& & & & & & \pncf & \!\mddots & & & \\
		\hdashline[2pt/2pt]
		& & & & & & \!\mddots &\pcf & \pcf & \mcdots & 0 \\		
		& & & & & \!\mddots & & & \pncf & & \\
		& & & & \!\mddots & & & & & \!\mddots & 1 \\		
    \end{array}
\right].
\end{equation}
The initial state $\pi_0$ and observation vector $u_\omega$ for a given $\omega$ are given by 
\begin{equation}
\pi_0 = \begin{bmatrix}
1 \\
0 \\
\mvdots \\
0
\end{bmatrix}, \quad
U^\prime_0 = \begin{bmatrix}
\overmat{\lambda^*}{1 & \mcdots & 1} & \overmat{(\lambda^*+1)l}{0 & \mcdots & 0}
\end{bmatrix} \quad \text{and}
\end{equation}
\vspace{0pt}
\begin{equation}
U^\prime_\omega = \begin{bmatrix}
\overmat{(\lambda^*+1)\omega-1}{0 & \mcdots & 0} & \overmat{\lambda^*+1}{1 & \mcdots & 1} & \overmat{(\lambda^*+1)(l-\omega)}{0 & \mcdots & 0}
\end{bmatrix} , \quad \omega = 1,\dots l.
\end{equation}

The dimension of the matrix $A$ scales with $M$ for both computations. This means that a straightforward algorithm for computing \eqref{eq:run_probability} will require $\mathcal{O}\left( M^3 \right)$ work. It is, however, only necessary to perform this computation once for each value of $\lambda^*\in\{1,\dots,\Lambda\}$ and $M$. Exploiting matrix sparseness and structure also reduces the total amount of computation. Given that we fix $M$ at level 1 for the entire simulation, it is possible to precompute a table with dimensions $\Lambda \times \frac{M}{2}$ containing all probabilities needed in the simulations. By tabulating the cumulative probabilities of the computed statistics in function of $\omega$, efficient inverse transform sampling can be implemented by taking a binary search approach. This can further be accelerated and made more cache efficient by starting the search at $\omega = \E \left[ \sigma_{\lambda^*} \right]$ and using linear interpolation to perform bisection steps in the search.

\section*{Acknowledgments}
We thank Ignace Bossuyt, Bert Mortier and Pieterjan Robbe for their valuable input on the topic of this work and proofreading efforts. We  also thank the anonymous reviewers for their detailed remarks on this manuscript.

\bibliographystyle{siamplain}
\bibliography{Combined}
\pagebreak
\includepdf[pages=-]{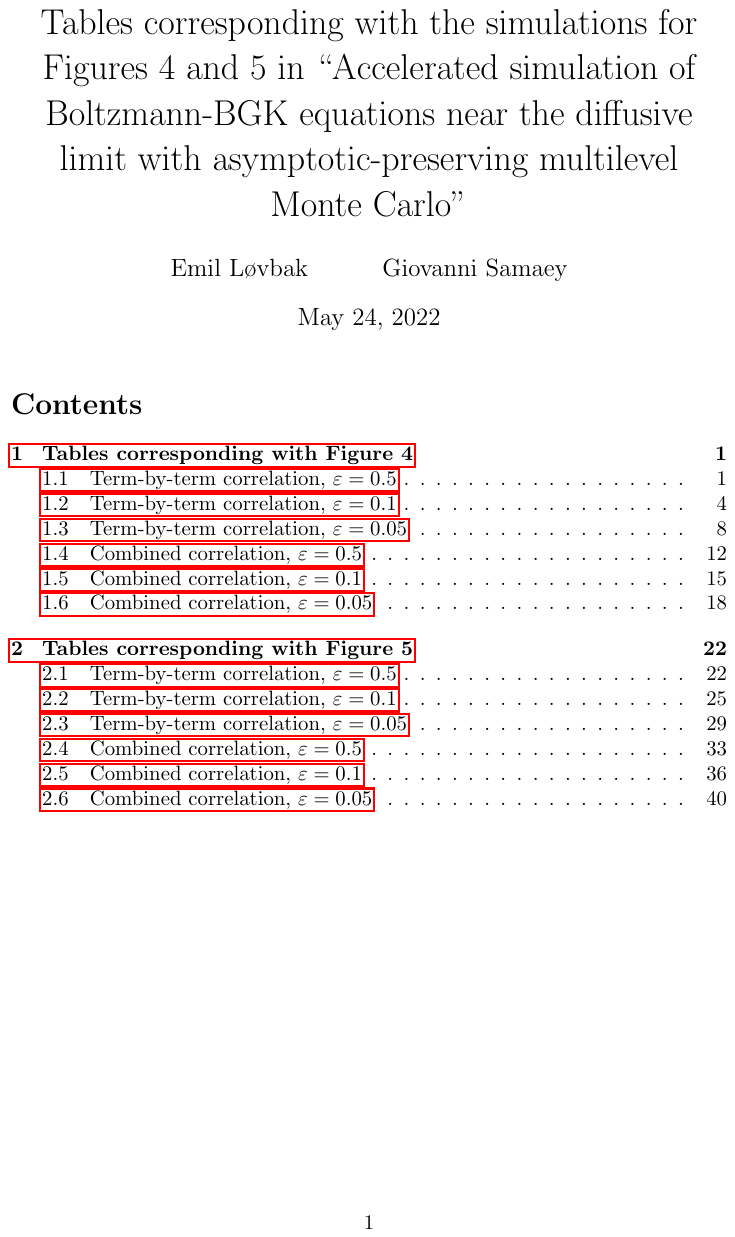}
\end{document}